\documentclass[11pt,english]{article}
\usepackage[T1]{fontenc}
\usepackage[latin9]{inputenc}
\usepackage{color}
\usepackage{float}
\usepackage{graphicx}
\usepackage{setspace}
\usepackage[authoryear]{natbib}
\usepackage[bookmarks=false]{hyperref}
\makeatletter

\usepackage{sectsty}
 
\usepackage{pslatex}

\sectionfont{\sffamily\large}
\subsectionfont{\sffamily\normalsize}
\newcounter{eqnum}[section] 
\newtheorem{proposition}{Proposition} 
 
 \newtheorem{theorem}{Theorem}
\newtheorem{lemma}{Lemma}

\newcommand{\eproof}{\mbox{}\hfill{\rule{8pt}{8pt}}}

\makeatother

\makeatother

\usepackage{babel}
\begin{document}

\title{Acquisition of Project-Specific Assets with Bayesian Updating\thanks{DOI: 10.1287/opre.1110.0949}}

\author{H. Dharma Kwon\thanks{Department of Business Administration, University of Illinois at Urbana-Champaign,
1206 South Sixth St., Champaign, IL 61820; Email: dhkwon@illinois.edu.}, Steven A. Lippman\thanks{UCLA Anderson School, UCLA, 110 Westwood Plaza, Los Angeles, CA 90095;
Email: slippman@anderson.ucla.edu}}

\date{September 12, 2010}
\maketitle
\begin{abstract}
We study the impact of learning on the optimal policy and the time-to-decision
in an\textcolor{red}{{} }\textcolor{black}{infinite-horizon }Bayesian
sequential decision model with two irreversible alternatives, exit
and expansion. In our model, a firm undertakes a small-scale pilot
project so as to learn, via Bayesian updating, about the project\textquoteright s
profitability, which is known to be in one of two possible states.
 The firm continuously observes the project\textquoteright s cumulative
profit, but the true state of the profitability is not immediately
revealed because of the inherent noise in the profit stream. The firm
bases its exit or expansion decision on the posterior probability
distribution of the profitability. The optimal policy is characterized
by a pair of thresholds for the posterior probability. We find that
the time-to-decision does not necessarily have a monotonic relation
with the arrival rate of new information. 
\end{abstract}
Time to decision. Brownian motion.\\
\emph{OR/MS subject classifications:} Decision analysis: sequential.
Finance: investment criteria. Probability: diffusion. 

\pagebreak

\section{Introduction}

Launching a new project entails enormous uncertainty. In order to
learn whether it is advisable to launch a new project, firms often
perform small-scale experiments before making an irreversible investment
in project-specific assets. One pertinent well-known example is US
steel maker Nucor's decision to adopt the world's first continuous
thin-slab casting technology. Recognizing the huge uncertainty in
the profitability of the new technology, Nucor built a pilot plant
in 1989 \citep{Ghemawat1998}. After the pilot plant proved to be
a success, Nucor expanded the use of the new technology by building
several other thin-slab casting plants beginning in 1992. In this
context of launching a new project, our paper explores the central
question of interest to the decision-maker: how uncertainty and learning
impact (1) the optimal expansion and exit policy and (2) the length
of time before an expansion or exit decision is made.

Our paper addresses three main difficulties with launching a new project
which requires project-specific investment. First, the prospect of
success is highly uncertain. Moreover, the firm undertaking the project
may have unproven capability. In the case of Nucor's adoption of the
thin-slab technology, there was uncertainty in the technical feasibility
as well as in the demand for the new product, and both factors contributed
to uncertainty in the profit. Second, when the uncertainty in the
profit is project-specific, it can be resolved only through experimentation.
Sources of uncertainty which are not project-specific can be resolved
by other means. Third, because the salvage value of the project-specific
assets can be very small, the expansion decision often is irreversible;
Nucor's adoption of thin-slab casting technology, which is useless
outside the industry, was irreversible. Hence, the decision to acquire
project-specific assets should be made only if the firm is sufficiently
confident of the project's success.

In this paper, we study a stylized model of a firm experimenting with
an unproven project. As the experimentation proceeds, the firm must
decide when to cease the experiment whether to expand or stop the
project. In our model, the profitability of the project is either
high or low, and it does not change in time. The firm knows the prior
(initial) probability that the new project's profitability is in the
high state. The firm experiments by launching a small scale enterprise.
Because the firm observes the profit at each point in time, it continuously
updates its belief regarding the profitability of the new project.
The firm cannot immediately determine the state of profitability because
the realized profit contains noise. At each point in time after the
launch, the firm can (1) continue the current experimentation, (2)
stop and exit the project, or (3) make an irreversible investment
to expand the project by acquiring project-specific assets.  

The cumulative profit is modeled as a Brownian motion with constant
drift and constant volatility $\sigma$. The drift of the Brownian
motion is known to be either in a high state or a low state, but the
true state of the drift is unknown. We formulate our model as a continuous-time
infinite-horizon optimal stopping problem. 

Employing the methods of optimal stopping, we obtain the following
results: (1) The optimal policy of expansion and exit is characterized
by two thresholds with respect to the posterior probability $P_{t}$
that the state is high: expand (abandon) the project if the posterior
probability $P_{t}$ hits the upper (lower) threshold; otherwise,
continue the pilot project. (2) The upper (lower) threshold decreases
(increases) in the volatility $\sigma$ of the cumulative profit.
(3) The expected time-to-decision has non-monotone dependence on $\sigma$. 

\citet{Shiryaev1967} and \citet{Peskir2006} studied the Bayesian
problem of minimizing the undiscounted cost of incorrect decisions
in sequential hypothesis testing in a setting in which the true state
is slowly revealed via a Brownian motion, and they obtained elegant
analytical solutions. In their model, the true state can be either
good or bad, and the decision to choose one hypothesis is irreversible.
They showed that the optimal stopping policy is characterized by a
pair of thresholds with respect to the posterior probability of being
in the good state. As explained later, they (and others) obtain a
stochastic differential equation for the posterior probability of
being in the good state at time $t$ given the history of the Brownian
motion up to time $t$. This stochastic differential equation is the
starting point of our analysis. Our model extends theirs by incorporating
discounting (i.e., the time value of money). In applying the model
to real options problems, it is essential to add discounting to the
model because models with discounting carry considerably more economic
interest and realism. In considering the generic issue of utilizing
a general purpose asset versus a specialized asset, we realized that
if the project turns out to be unprofitable on an operating basis,
the purchaser of the specialized asset would have bought a white elephant,
an asset which has no economic use. However, even if the project turns
out to be unprofitable, the purchaser might have the option to exit
the project, albeit having sunk funds into the purchase of a specialized
asset which no longer has any economic value. For this reason, we
also modify Shiryaev's model by incorporating an embedded exit option
after expansion (see Sec. \ref{sec:Exit-Option}). Our principal goals
are to study the comparative statics of the time-to-decision in our
Bayesian real options model as well as to characterize the optimal
policy. 

Although ours is an extension of Shiryaev's model, our model addresses
a different economic problem, and its comparative statics analysis
is different from that of Shiryaev's model. Our model does not have
a zero-discount rate counterpart because it allows perpetual operation
of a project after expansion whence it does not simply reduce to Shiryaev's
model as the discount rate goes to zero. Also, due to the difference
in the economic context, Shiryaev's model does not allow for an embedded
option to stop after expansion as our model does. Furthermore, the
comparative statics analysis of our model with respect to the volatility
is more complex  because of the non-zero discount rate which adds
an additional model parameter, and the large-volatility results are
different because of the difference in the economic problem. For example,
the expected time-to-decision is never monotonically increasing in
the volatility in Shiryaev's model whereas there are parameter values
in our model under which the relationship is monotone (see Figure
2).

Our paper contributes to the literature on real options theory, in
particular, in the context of the comparative statics with respect
to the volatility. The theory of real options has shown that there
is value to waiting before making an irreversible decision when the
value of such a decision is uncertain. Moreover, both the asset value
and the value of waiting increase in the volatility (\citealt{Dixit1992,Alvarez2003}).
In our model, the volatility of $P_{t}$ decreases in $\sigma$: as
$\sigma$ increases, the arrival of new information slows down as
shown by \citet{Bergemann2000}. The asset value, that is, the value
of the project and the associated option to expand or exit, decreases
in $\sigma$ because the arrival of information about the true state
slows down. Likewise, the value of waiting also decreases in $\sigma$.
A decrease in the value of waiting also decreases the time-to-decision.
There is, however, a countervailing effect of an increase in $\sigma$.
In our model, $\sigma$ constitutes noise in the observed profit stream,
so the arrival rate of information is higher when $\sigma$ is low.
If $\sigma$ increases, then the decision-maker might want to wait
a very long time to collect enough information to increase the likelihood
of making the correct decision. To our knowledge, our paper is the
first work that studies the resulting effect of an increase in $\sigma$
on the time-to-decision, which depends on the relative magnitudes
of these two countervailing effects. 

We study the dependence of the expected time-to-decision on $\sigma$
and obtain the following results. For sufficiently small values of
$\sigma$, the effect of the rate of information arrival dominates
the effect of the value of waiting, so the expected time-to-decision
increases in $\sigma$. In contrast, for sufficiently large values
of $\sigma$, the resulting comparative statics of the expected time-to-decision
with respect to $\sigma$ depends on the model parameters. 

Our results also have practical implications for firms faced with
expansion and exit decisions. In most business decisions, the project's
lifetime is a significant factor. However, there is a paucity of work
on the length of time-to-decision; our paper begins to fill this gap
in the literature. 

The paper is organized as follows. We review related literature in
Sec. \ref{sec:Literature}. We solve and analyze the infinite-horizon
model of a firm learning from a pilot project with an exit and expansion
decision in Sec. \ref{sec:The-Model}. In Sec. \ref{sec:Exit-Option},
we allow the decision-maker to exit even after the expansion, and
we study the effect of the post-expansion exit option. The implications
of having an additional source of information are discussed in Sec.
\ref{sec:Discussion}. Conclusions are given in Sec. \ref{sec:Conclusions}. 

\section{Related Literature \label{sec:Literature}}

Our paper is closely related to the literature on the value of acquiring
information. The seminal work on the subject of learning-by-doing
is \citet{Arrow1962}. In particular, there is a strand of papers
that analyzed optimal stopping models to study the value of experimentation
under incomplete information. \citet{Bernardo2002} used numerical
methods to study the real option value of a firm which invests in
either a specialized or generalized project, learns about its own
capabilities, and then decides either to scale up the specialized
project or to expand into a multisegment business. \citet{McCardle1985}
and \citet{Ulu2009} studied the adoption and abandonment  of new
technologies under incomplete information regarding profits when the
information acquisition cost per unit time is deterministic; they
characterized the optimal adoption policy, but they did not study
the time-to-decision. In our model, the implicit cost of acquiring
information is uncertain.

Our paper is an extension of the classic sequential hypothesis testing
problem, the objective of which is to minimize the expected cost of
errors (cost of choosing the incorrect hypothesis) when the probability
distribution of the hypotheses is updated in a Bayesian manner. The
subject of sequential hypothesis testing was pioneered by \citet{Wald1945}
whose work spawned a vast literature on this subject. See \citet{Poor2008}
and \citet{Lai2001} for references. 

Sequential hypothesis testing problems can be solved as stopping time
problems. There is a vast literature on optimal stopping problems,
and the standard exposition of the optimal stopping theory in continuous
time can be found in \citet{Shiryaev1978}, \citet{Oksendal03}, and
\citet{Peskir2006}. The solution method consists of solving a characteristic
partial differential equation and applying smooth-pasting conditions;
the existence of such a solution automatically ensures an optimal
solution (\citealt{Dayanik2003}, Chapter 10 of \citealt{Oksendal03},
and Chapter IV of \citealt{Peskir2006}). (The continuous-time optimal
stopping theory has been applied to study mathematical properties
of real options models; see, for example, \citealt{Alvarez2001} and
\citealt{Wang2005}.) We solve our model by formulating it as an optimal
stopping problem. In solving it, we build upon the work of \citet{Shiryaev1967}
and \citet{Peskir2006}. \citet{Shiryaev1967} studied the problem
of minimizing the cost of error with two simple hypotheses on the
drift of a one-dimensional Brownian motion, and he obtained an analytical
solution. The basic element of his model is the stochastic differential
equation for the posterior probability distribution. This is the starting
point of our analysis. 

This two-drift Brownian motion model has been applied in many economics
papers concerning the optimal level of experimentation in several
different contexts. (See \citealt{Keller1999}, \citealt{Moscarini2001},
\citealt{Bolton1999}, and \citealt{Bergemann2000}.) For the most
part, these papers focus upon decision-making which leads to more
rapid learning.  In contrast, our paper focuses upon the optimal stopping
decision.

In an optimal stopping problem which applied Shiryaev's framework,
\citet{Ryan2003} considered when to stop (abandon) a project with
unknown profitability. In their problem, a firm seeks to maximize
the discounted cumulative profit from a project. The firm can abandon
the project at any point in time, and the cumulative profit is modeled
as a Brownian motion in which the drift takes either a known positive
value or a known negative value. Employing the methods of stochastic
analysis, they obtained an optimal policy which is stationary and
characterized by a threshold on the posterior probability.  \citet{Decamps2005}
also employed Shiryaev's framework to study the optimal time to invest
in an asset with an unknown underlying value. In their model, the
reward from stopping is the Brownian motion itself rather than the
expected value of the time-integral of a Brownian motion. In both
of these papers, there is a single alternative to continuing, and
the optimal policies are characterized by a single threshold. 

In our model, there are two alternatives to continuing, and the optimal
policy entails two thresholds. This result is a similar to the two-threshold
policy obtained by \citet{Shiryaev1967} and \citet{Peskir2006}.
Policies with two thresholds also occur in many models under \emph{complete}
information. See, for example, \citet{Alvarez2003a} and \citet{Decamps2006}.

Finally, the comparative statics properties of optimal stopping policies
with respect to the volatility are of particular interest. In the
economics literature, it has been shown that the value of waiting
before decision has a monotonicity property with respect to the volatility
(\citealt{Dixit1992}). In an economics model of a firm which has
an option to enter and exit an industry, \citet{Dixit1989} employed
numerical examples to depict the comparative statics of the optimal
entry and exit thresholds with respect to the volatility of the profit
stream. In their Bayesian decision model, \citet{Bernardo2002} demonstrated
by numerical illustration that the volatility of the observed profit
shrinks the continuation region. \citet{Alvarez2003} obtained general
comparative statics result for the optimal policy and the optimal
return with respect to the volatility for a class of optimal stopping
problems which arise in economic decisions. \citet{Kwon2010} showed
that the comparative statics of the optimal policy with respect to
the volatility is non-trivial when there is an embedded option. All
of these papers focus on the effect of the volatility on the optimal
policy or the optimal return; none of these papers studied the impact
of volatility on the time-to-decision. In the classical sequential
hypothesis testing problem, \citeauthor{Wald1973} (1973, Chapter
3) studied the expected number of observations before choosing a hypothesis,
but its dependence on the noise level of the observation was not explored.

\section{Base Model\label{sec:The-Model}}

Consider the decision problem of a firm launching a new project. The
state $S$ of the project lies in $\{H,L\}$: the project will bring
in either a high profit stream if it is in state $H$ or a low profit
stream if it is in state $L$.  The state is not known, but the firm
knows $p_{0}$, the prior probability that the state is $H$. Before
launching the project on a large scale, the firm can collect information
about the profitability of the project by operating a pilot project.
During the pilot phase, the mean profit per unit time is either $h$
if $S=H$ or $\ell$ if $S=L$. We assume that $h>0$ and $\ell<0$
so the project is profitable only if it is in the high state. At each
point in time during the pilot phase, the firm can (1) expand the
project by acquiring project-specific assets (equipment) at cost $k$
or (2) permanently exit the project. Once the project is expanded,
the mean profit stream is changed by $h^{\dagger}$ in state $H$
and $\ell^{\dagger}$ in state $L$ so that the mean profit rate after
acquisition is $h^{A}\equiv h+h^{\dagger}$ or $\ell^{A}\equiv\ell+\ell^{\dagger}$.
Let $\alpha>0$ denote the continuous time discount rate.  In this
section, we limit consideration to the case in which exit is never
optimal after expansion due to a prohibitively high exit cost. In
similar settings with investment in project-specific assets, the assumption
of a prohibitively high exit cost was also made by \citet{Lippman1992}
and \citet{Bernardo2002}. We relax this assumption in Sec. \ref{sec:Exit-Option}
and investigate its effect. 

The firm observes the profit stream without error. However, because
of noise, the firm can never perfectly determine the underlying profit
rate (i.e., state of the project). We model the cumulative profit
stream $\{X_{t}:t\ge0\}$ during the pilot phase as a Brownian motion
with  unknown drift:
\begin{equation}
X_{t}=\mu t+\sigma B_{t}\:,\label{eq:X-t}
\end{equation}
where the unknown drift $\mu\in\{h,\ell\}$ is the mean profit per
unit time, $\sigma$ is the volatility of the cumulative profit, and
$\mathbf{B}=\{B_{t}:t\ge0\}$ is a one-dimensional standard Brownian
motion. After the firm expands the project, the new profit stream
$\{Y_{t}:t\ge0\}$ is given by $Y_{t}=\mu^{A}t+\sigma B_{t}$ where
$\mu^{A}\in\{h^{A},\ell^{A}\}$; if the firm exits, then the profit
stream is zero. 

We use a continuous-time approach and model the profit stream as a
one-dimensional Brownian motion because of the analytical tractability
of the resulting model. Likewise, numerous papers in economics have
 modeled cumulative profits under incomplete information as Brownian
motion with unknown drift and obtained useful insights from analytically
tractable models (\citealt{Bolton1999,Keller1999,Bergemann2000,Bernardo2002,Ryan2003,Decamps2005}).

This model is general enough to be applicable in many contexts.

\noindent \textbf{Example 1 (Expansion):} Consider a firm testing
a new unproven technology by operating a single pilot production plant.
If the plant proves to be a success, $n$ more plants will be built
at cost $k$. Alternatively, $n$ can be regarded as an increase in
the production capacity. Suppose that the profit per plant is not
changed by the expansion. The expected profit per unit time from the
pilot plant is $h>0$ or $\ell<0$, and the expected profit per unit
time from the $n+1$ plants after expansion will be $h^{A}=(n+1)h$
or $\ell^{A}=(n+1)\ell$. \eproof

\noindent \textbf{Example 2 (Acquisition of a dedicated facility):}
Consider a firm testing a new project by outsourcing the main task
or by renting a facility for the core task of the project. By outsourcing
or renting, the firm can easily terminate the project, but the profit
stream is diminished because it has to pay a price for outsourcing
or renting. The firm can acquire a dedicated facility of its own to
avoid paying the price for renting or outsourcing. Before the acquisition
of the dedicated facility, the expected profit per unit time is $h$
or $\ell$; after the acquisition, the expected profit per unit time
improves by $h^{\dagger}$ and $\ell^{\dagger}$ which are the expected
cost of outsourcing or renting in the high and low states respectively.
The increase $h^{\dagger}$ or $\ell^{\dagger}$ in expected profit
per unit time might be state-dependent because the cost of renting
might be on a per unit sold basis. \eproof 

To simplify our presentation, we will focus on expansion decisions
(Example 1) with the understanding that our model fully addresses
Example 2 as well. For consistency, we shall speak in terms of project
expansion rather than acquisition of project-specific assets.

\subsection{Posterior Probability\label{subsec:posterior}}

At time zero, the firm's prior probability that $S=H$ is $p_{0}$.
Our first task is to compute the posterior probability that the project
is in the high state at time $t$. Subsequently, we utilize the posterior
probability process to ascertain the optimal expansion and exit decisions. 

Let $X_{t}$, $\mu$, and $B_{t}$ of Eq. (\ref{eq:X-t}) be defined
on a probability space $(\Omega,\mathcal{G},\mathcal{P})$ so that
$X_{t}$, $\mu$, and $B_{t}$ are measurable with respect to $\mathcal{G}$,
and let $\{\mathcal{F}_{t}:t\ge0\}$ be the natural filtration with
respect to the observable process $\{X_{t}:t\ge0\}$. (The unknown
drift $\mu$ and the unobservable process $\{B_{t}:t\ge0\}$ are \emph{not}
adapted to the filtration $\mathcal{F}_{t}$.) Here the probability
measure is denoted by $\mathcal{P}$. Let $P_{t}\equiv\mathcal{P}[\mu=h\vert\mathcal{F}_{t}]=E^{p_{0}}[1_{\{\mu=h\}}\vert{\mathcal{F}}_{t}]$
denote the posterior probability of the event $\{\mu=h\}$ at time
$t$, where $E^{p_{0}}[\cdot]$ denotes the expectation conditional
on the initial condition $P_{0}=p_{0}$, and $1_{A}$ is the indicator
function of a set $A$. The posterior probability process ${\bf P}$
is a continuous martingale with respect to $\{{\mathcal{F}}_{t}:t\ge0\}$,
a strong Markov process that is homogeneous in time, and the unique
solution to the following stochastic differential equation (Liptser
and Shiryaev 1974, p. 371; Peskir and Shiryaev 2006, pp. 288-289;
Ryan and Lippman 2003, pp. 252-254): 
\begin{equation}
dP_{t}=\frac{h-\ell}{\sigma}P_{t}(1-P_{t})d\tilde{B}_{t},\label{eq:P-t-SDE}
\end{equation}
 where $\tilde{{\bf B}}=\{\tilde{B}_{t},{\mathcal{F}}_{t}:t\ge0\}$
is a one-dimensional standard Brownian motion defined by 
\[
\tilde{B}_{t}\equiv\frac{1}{\sigma}\left(X_{t}-\int_{0}^{t}E[\mu\vert{\mathcal{F}}_{s}]ds\right)=\frac{1}{\sigma}[X_{t}-\int_{0}^{t}(P_{s}h+(1-P_{s})\ell)ds]\:.
\]
Unlike $B_{t}$ which is unobservable, the new Brownian motion $\tilde{B}_{t}$
is observable because it can be completely constructed from the observed
value of $X_{t}$ over time.

Equation (\ref{eq:P-t-SDE}) implies that the posterior process undergoes
a more rapid change if the coefficient $\frac{h-\ell}{\sigma}P_{t}(1-P_{t})$
increases. The factor $\frac{h-\ell}{\sigma}$, which is the uncertainty
in the drift normalized by the volatility, has the interpretation
of the signal-to-noise ratio. Hence, $\frac{h-\ell}{\sigma}$ measures
the rate of arrival of new information, and it follows that the posterior
evolves more rapidly with a high value of $\frac{h-\ell}{\sigma}$
(\citealt{Bergemann2000}). The additional factor $P_{t}(1-P_{t})$,
which peaks at $P_{t}=0.5$, results from the fact that the posterior
process requires more information (and hence a longer time) to change
if $P_{t}$ is closer to strong beliefs (either $0$ or $1$). 

Using the strong Markov property of the process ${\bf P}$ and Bayes
rule (Peskir and Shiryaev, 2006, pp. 288-289), we have 
\begin{eqnarray}
P_{t} & \equiv & E^{p_{0}}[1_{\{\mu=h\}}\vert{\mathcal{F}}_{t}]=E^{p_{0}}[1_{\{\mu=h\}}\vert X_{t}]=\frac{p_{0}\exp\{-\frac{(X_{t}-ht)^{2}}{2\sigma^{2}t}\}}{p_{0}\exp\{-\frac{(X_{t}-ht)^{2}}{2\sigma^{2}t}\}+(1-p_{0})\exp\{-\frac{(X_{t}-\ell t)^{2}}{2\sigma^{2}t}\}}\nonumber \\
 & = & \left(1+\frac{1-p_{0}}{p_{0}}\exp\{-\frac{h-\ell}{\sigma^{2}}\cdot[(\mu-\frac{h+\ell}{2})t+\sigma B_{t}]\}\right)^{-1}\;.\label{eq:Pt-formula}
\end{eqnarray}
Because $\mu-\frac{h+\ell}{2}=(h-\ell)/2$ if $\mu=h$ and $\mu-\frac{h+\ell}{2}=(\ell-h)/2$
if $\mu=\ell$, Eq. (\ref{eq:Pt-formula}) reflects the fact that
$P_{t}$ tends to increase if $\mu=h$ and decrease if $\mu=\ell$. 

\subsection{The Objective Function \label{subsec:Objective}}

The firm seeks to maximize its expected discounted cumulative profit
where $\alpha>0$ is the discount rate. Suppose the firm decides to
expand at the stopping time $\tau$. Then the discounted reward over
the time interval $(\tau,\infty)$ discounted back to time $\tau$
is given by 
\[
E^{p_{0}}[\int_{\tau}^{\infty}e^{-\alpha(t-\tau)}dY_{t}\vert\mathcal{F}_{\tau}]=E^{P_{\tau}}[\int_{0}^{\infty}e^{-\alpha t}dY_{t}]=E^{P_{\tau}}[\int_{0}^{\infty}\mu^{A}e^{-\alpha t}dt]=\frac{1}{\alpha}[P_{\tau}h^{A}+(1-P_{\tau})\ell^{A}]\;,
\]
where the first equality follows from  $\mathbf{P}$ being a strong
Markov process and the second from $E^{x}[\int_{0}^{\infty}e^{-\alpha t}dB_{t}]=0$.
If the firm decides to exit, then the associated reward from exit
is zero. Thus, the optimal reward $r(\cdot)$ from stopping at time
$\tau$ is the greater of the reward from expansion and the reward
from exit:
\[
r(P_{\tau})=\max\{\frac{1}{\alpha}[P_{\tau}h^{A}+(1-P_{\tau})\ell^{A}]-k,0\}\:,
\]
where $k$ is the cost of expansion. 

The objective function of the firm is the time-integral of the discounted
profit stream up to time $\tau$ plus the reward from stopping at
time $\tau$:
\begin{eqnarray}
V_{\tau}(p_{0}) & = & E^{p_{0}}[\int_{0}^{\tau}e^{-\alpha t}dX_{t}+e^{-\alpha\tau}r(P_{\tau})]=E^{p_{0}}[\frac{\mu}{\alpha}(1-e^{-\alpha\tau})+e^{-\alpha\tau}r(P_{\tau})]\nonumber \\
 & = & \frac{1}{\alpha}[p_{0}h+(1-p_{0})\ell]+E^{p_{0}}[e^{-\alpha\tau}g(P_{\tau})]\:,\label{eq:V}
\end{eqnarray}
where 
\begin{equation}
g(x)=\frac{1}{\alpha}\max\{xh^{\dagger}+(1-x)\ell^{\dagger}-k\alpha,-xh-(1-x)\ell\}\:.\label{eq:g-x}
\end{equation}
The reward function $g(x)=V_{0}(x)-V_{\infty}(x)$ is the difference
between the return from immediate stopping ($\tau=0$) and the return
from never stopping ($\tau=\infty$).

Because the only $\tau$-dependence of $V_{\tau}(p_{0})$ is in the
term $E^{p_{0}}[e^{-\alpha\tau}g(P_{\tau})]$, the stopping problem
of Eq. (\ref{eq:V}) is equivalent to maximizing 
\begin{equation}
R_{\tau}(p_{0})=E^{p_{0}}[e^{-\alpha\tau}g(P_{\tau})]\:.\label{eq:R-p}
\end{equation}
For convenience, we regard $R_{\tau}(p_{0})$ as the objective function
for the remainder of this section.

\subsection{Optimal Policy \label{subsec:Policy}}

As revealed in Eq. (\ref{eq:V}), the firm seeks a stopping time $\tau^{*}$,
termed optimal, that maximizes $R_{\tau}(p_{0})$ defined in Eq. (\ref{eq:R-p}).
Because our infinite horizon stopping problem is stationary, it comes
as no surprise that it suffices to focus on the class of stopping
times that are time-invariant (\citealt{Oksendal03}, p. 220). Each
time-invariant stopping time can be defined through a set $A$ as
follows:
\[
\tau_{A}\equiv\inf\{t\ge0:P_{t}\not\in A\}\:.
\]
We call $\tau_{A}$ \emph{the exit time from the set} $A$. For notational
convenience, we let $R^{*}(p)\equiv R_{\tau^{*}}(p)$ denote the optimal
return function whenever the optimal stopping time $\tau^{*}$ exists.

The existence of an optimal stopping time is not guaranteed in general,
but we can prove that our model satisfies the sufficient conditions
for the existence of an optimal stopping time. Before proving the
sufficient conditions, we need to lay out some technical preliminaries.
Clearly, if $P_{t}=p$ at some time $t$, then it is optimal to continue
if the return $R^{*}(p)$ from continuing exceeds the return $R_{0}(p)=g(p)$
from stopping. Accordingly, we refer to $\{p:R^{*}(p)>R_{0}(p)\}$
as the \emph{continuation region}. Let 
\[
\mathcal{A}\equiv-\alpha+\frac{1}{2}(\frac{h-\ell}{\sigma})^{2}p^{2}(1-p)^{2}\partial_{p}^{2}
\]
be the characteristic differential operator for $P_{t}$; $\mathcal{A}$
plays a fundamental role in the solution of stopping time problems.
Here $\partial_{t}$ is replaced by $-\alpha$ because our solution
is time-invariant except for the discount factor $e^{-\alpha t}$.
One of the sufficient conditions for the optimality of $\tau^{*}$
is that $R^{*}(\cdot)$ satisfy the differential equation $\mathcal{A}R^{*}(p)=0$
if $p$ belongs to the continuation region. There are two \emph{fundamental
solutions,} $\psi$ and $\phi$, to the second-order linear ordinary
differential equation $\mathcal{A}f(p)=0$: 
\begin{eqnarray*}
\psi(p) & = & p^{\frac{1}{2}(1+\gamma)}(1-p)^{\frac{1}{2}(1-\gamma)}\:,\\
\phi(p) & = & p^{\frac{1}{2}(1-\gamma)}(1-p)^{\frac{1}{2}(1+\gamma)}\:,
\end{eqnarray*}
where 
\begin{equation}
\gamma\equiv\sqrt{1+\frac{8\alpha\sigma^{2}}{(h-\ell)^{2}}}\:.\label{eq:gamma}
\end{equation}
Note that $\psi(\cdot)$ is convex increasing and $\phi(\cdot)$ is
convex decreasing.

If the firm stops at time $t$ with $P_{t}=p$, then the firm's (expected
discounted) return over the interval $[t,\infty)$ is $g(p)$. From
Eq. (\ref{eq:g-x}), note that $g(p)$ is the maximum of two functions
$g_{1}(p)\equiv ph^{\dagger}/\alpha+(1-p)\ell^{\dagger}/\alpha-k$
and $g_{2}(p)\equiv-ph/\alpha-(1-p)\ell/\alpha$, each of which is
linear in the posterior probability $P_{\tau^{*}}=p$ at the stopping
time. The first function $g_{1}(p)$ gives the return associated with
expansion and is increasing in $p$. The second function $g_{2}(p)$
is associated with the return to exit and is decreasing in $p$. if
$\frac{k\alpha-\ell^{A}}{h^{A}-\ell^{A}}\in(0,1)$, then the two linear
functions cross at $\hat{p}$ where 
\[
\hat{p}\equiv\frac{k\alpha-\ell^{A}}{h^{A}-\ell^{A}}\:,
\]
and $g(\cdot)$ is minimized at $\hat{p}$. Note that $g_{1}(p)>g_{2}(p)$
if and only if $p>\hat{p}$, so it is optimal to expand if $P_{\tau^{*}}\ge\hat{p}$
and optimal to exit if $P_{\tau^{*}}<\hat{p}$. 

First, we consider the cases $\frac{k\alpha-\ell^{A}}{h^{A}-\ell^{A}}\ge1$
or $\frac{k\alpha-\ell^{A}}{h^{A}-\ell^{A}}\le0$.

\begin{proposition} \label{prop:singe-threshold} (i) Suppose $\frac{k\alpha-\ell^{A}}{h^{A}-\ell^{A}}\ge1$.
Then the optimal policy is to exit when $P_{t}$ hits a lower threshold
\[
\underline{\theta}=\frac{-(\gamma-1)\ell}{(\gamma+1)h-(\gamma-1)\ell}\:
\]
 and to continue operation otherwise. The optimal return is given
by
\[
R^{*}(p)=\left\{ \begin{array}{ll}
\frac{-2\ell h/\alpha}{\sqrt{(\gamma^{2}-1)(-\ell h)}}[\frac{(\gamma-1)}{(\gamma+1)}(\frac{-\ell}{h})]^{\frac{\gamma}{2}}\cdot\phi(p) & \mbox{if }\:p>\underline{\theta}\\
g(p) & \mbox{otherwise}\:.
\end{array}\right.
\]
(ii) Suppose $\frac{k\alpha-\ell^{A}}{h^{A}-\ell^{A}}\le0$. Then
the optimal policy is to expand the project when $P_{t}$ hits an
upper threshold 
\[
\bar{\theta}=\frac{-(\gamma+1)(\ell^{\dagger}-k\alpha)}{(\gamma-1)(h^{\dagger}-k\alpha)-(\gamma+1)(\ell^{\dagger}-k\alpha)}\:
\]
 and to continue operation otherwise. The optimal return is given
by
\[
R^{*}(p)=\left\{ \begin{array}{ll}
\frac{-2(\ell^{\dagger}-k\alpha)(h^{\dagger}-k\alpha)/\alpha}{\sqrt{(\gamma^{2}-1)(h^{\dagger}-k\alpha)(k\alpha-\ell^{\dagger})}}[\frac{(\gamma-1)}{(\gamma+1)}(\frac{h^{\dagger}-k\alpha}{k\alpha-\ell^{\dagger}})]^{\frac{\gamma}{2}}\cdot\psi(p) & \mbox{if }\:p<\bar{\theta}\\
g(p) & \mbox{otherwise}\:.
\end{array}\right.
\]
\end{proposition}

If $\frac{k\alpha-\ell^{A}}{h^{A}-\ell^{A}}\ge1$ or $k\alpha\ge h^{A}$,
the exit option is always better than the expansion option because
the expansion cost is too high, so expansion is never optimal. If
$\frac{k\alpha-\ell^{A}}{h^{A}-\ell^{A}}\le0$ or $k\alpha\le\ell^{A}$,
then the expansion option is always better than the exit option because
the profit improvement, even in the low state, exceeds the cost of
investment. In these two cases, the optimal policy is characterized
by only one threshold. The single-threshold solutions are studied
in detail by \citet{Ryan2003} and \citet{Kwon2008a}. Henceforth,
we focus on the more interesting case of $\hat{p}=\frac{k\alpha-\ell^{A}}{h^{A}-\ell^{A}}\in(0,1)$
when the optimal policy is characterized by \emph{two} thresholds. 

\begin{theorem} \label{thm:basic-policy} Assume $\frac{k\alpha-\ell^{A}}{h^{A}-\ell^{A}}\in(0,1)$.
(i) The optimal policy to Eq. (\ref{eq:R-p}) always exists, and there
is a pair of thresholds $\underline{\theta}$ and $\bar{\theta}$
which satisfy $0<\underline{\theta}<\bar{\theta}<1$ such that the
optimal stopping time is given by 
\[
\tau^{*}=\inf\{t>0:P_{t}\not\in(\underline{\theta},\bar{\theta})\}\:.
\]
 The optimal return function is given by
\begin{equation}
R^{*}(p)=\left\{ \begin{array}{ll}
c_{1}\psi(p)+c_{2}\phi(p) & \mbox{if }\:p\in(\underline{\theta},\bar{\theta})\\
g(p) & \mbox{otherwise}\:,
\end{array}\right.\label{eq:Optimal-return}
\end{equation}
 where $c_{1}$ and $c_{2}$ are the unique positive numbers such
that $R^{*}(\cdot)$ is continuously differentiable. 

(ii) The time-to-decision $\tau^{*}$ is finite with probability one,
and 
\begin{eqnarray}
E^{p}[\tau^{*}] & = & \frac{p-\underline{\theta}}{\overline{\theta}-\underline{\theta}}T(\bar{\theta})+\frac{\bar{\theta}-p}{\overline{\theta}-\underline{\theta}}T(\underline{\theta})-T(p)\label{eq:E-tau}
\end{eqnarray}
where $T(p)\equiv\frac{2\sigma^{2}}{(h-\ell)^{2}}(2p-1)\ln(\frac{p}{1-p})$.
\end{theorem}

The proof of Theorem \ref{thm:basic-policy} is given in the Appendix.
The proof of the existence of an optimal policy consists of constructing
a function $f(\cdot)$ (a candidate for the optimal return function)
which satisfies the differential equation $\mathcal{A}f(p)=0$, the
boundary conditions $f(p)=g(p)$ for $p\in\{\underline{\theta},\bar{\theta}\}$,
and continuous differentiability (smooth pasting conditions). The
expected time-to-decision formula in Eq. (\ref{eq:E-tau}) was also
proven in \citet{Poor2008}, p. 87,  using the standard property of
martingales; we reproduce its proof in the Appendix.

The optimal policy is intuitively straightforward: expand the project
if the posterior probability $P_{t}$ is high enough, and exit if
it is low enough; otherwise continue the pilot project. The expected
time-to-decision $E^{p}[\tau^{*}]$ is finite because the optimal
policy is characterized by two thresholds, either of which is reached
by $P_{t}$ eventually. By contrast, in exit-only or investment-only
models as in the case of Proposition \ref{prop:singe-threshold} or
as in \citet{Ryan2003}, the expected time to decision is infinite
as there is a positive probability that the threshold is never reached
in a finite amount of time. 

The optimal return $R^{*}(\cdot)$ given in Eq. (\ref{eq:Optimal-return})
is strictly convex and larger than $g(p)$ for $p\in(\underline{\theta},\bar{\theta})$
which is clear from the fact that $\psi$ and $\phi$ are convex and
$c_{1}$ and $c_{2}$ are positive. 

\subsection{Comparative Statics \label{subsec:ComparativeStatics}}

The analytical solution obtained in the previous subsection enables
us to obtain comparative statics with respect to $\sigma$. 

\begin{proposition} \label{prop:sigma} The optimal return function
$V_{\tau^{*}}(p)$ is non-increasing in $\sigma$, and the upper (lower)
threshold decreases (increases) in $\sigma$. \end{proposition}

\noindent 

As the noise level $\sigma$ of the observed profit increases, the
arrival of new information slows down, and it takes longer to accumulate
information about the profit state. It follows that it is optimal
to set a smaller continuation region (smaller values of the upper
threshold and larger values of the lower threshold) as $\sigma$ increases.
Thus, the upper threshold decreases and the lower threshold increases
in $\sigma$. 

Proposition \ref{prop:sigma} is entirely consistent with the results
from conventional real options theory. It is well-known that the value
of a real option increases in the volatility of the asset value process
and that the continuation region grows with the volatility. In our
Bayesian stopping problem, by Eq. (\ref{eq:P-t-SDE}), $\frac{h-\ell}{\sigma}P_{t}(1-P_{t})$
is the volatility of the posterior process $\{P_{t}:t\ge0\}$. As
per Proposition \ref{prop:sigma}, an increase in the volatility $\frac{h-\ell}{\sigma}P_{t}(1-P_{t})$
(i.e., an increase in the rate of information arrival) leads to an
increase in the value of the optimal return $V_{\tau^{*}}(\cdot)$.

Because we do not have closed-form expressions for $\underline{\theta}$
and $\bar{\theta}$, Eq. (\ref{eq:E-tau}) does not provide a closed-form
expression for the expected time-to-decision as a function of $\sigma$.
However, asymptotic results can be derived from Eq. (\ref{eq:E-tau}).
Theorem \ref{thm:Asymptotics} provides the asymptotic properties
of the thresholds as well as the expected time to decision for small
and large values of $\sigma$.

\begin{theorem}[Expected Time-to-Decision] \label{thm:Asymptotics}
(i) As $\sigma\rightarrow0$, $\underline{\theta}\searrow0$, $\bar{\theta}\nearrow1$,
and $E^{p}[\tau^{*}]\searrow0$. (ii) Suppose $g(\hat{p})\ge0$. As
$\sigma\rightarrow\infty$, $\underline{\theta}\nearrow\hat{p}$,
$\bar{\theta}\searrow\hat{p}$, and $E^{\hat{p}}[\tau^{*}]\searrow0$.
(iii) Suppose $g(\hat{p})<0$. As $\sigma\rightarrow\infty$, $\underline{\theta}\nearrow\frac{-\ell}{h-\ell}$,
$\bar{\theta}\searrow\frac{-\ell^{\dagger}+k\alpha}{h^{\dagger}-\ell^{\dagger}}$,
$E^{p}[\tau^{*}]\rightarrow\infty$ if $p\in(\frac{-\ell}{h-\ell},\frac{-\ell^{\dagger}+k\alpha}{h^{\dagger}-\ell^{\dagger}})$,
$E^{p}[\tau^{*}]\rightarrow0$ if $p\not\in[\frac{-\ell}{h-\ell},\frac{-\ell^{\dagger}+k\alpha}{h^{\dagger}-\ell^{\dagger}}]$,
and $E^{p}[\tau^{*}]=O(\sigma)$ if $p\in\{\frac{-\ell}{h-\ell},\frac{-\ell^{\dagger}+k\alpha}{h^{\dagger}-\ell^{\dagger}}\}$.\footnote{A function $f(\sigma)$ is said to be $O(\sigma)$ if there is a positive
constant $M$ such that $\vert f(\sigma)\vert<M\sigma$ for all $\sigma$
sufficiently large.}\end{theorem}

Note that $g(p)<0$ if and only if $\frac{-\ell}{h-\ell}<p<\frac{-\ell^{\dagger}+k\alpha}{h^{\dagger}-\ell^{\dagger}}$.
If $g(\hat{p})<0$, the interval $(\frac{-\ell}{h-\ell},\frac{-\ell^{\dagger}+k\alpha}{h^{\dagger}-\ell^{\dagger}})$
always remains a subset of the continuation region so $\underline{\theta}<\frac{-\ell}{h-\ell}$
and $\frac{-\ell^{\dagger}+k\alpha}{h^{\dagger}-\ell^{\dagger}}<\bar{\theta}$
because an immediate stop in the interval $(\frac{-\ell}{h-\ell},\frac{-\ell^{\dagger}+k\alpha}{h^{\dagger}-\ell^{\dagger}})$
would result in a negative reward.

\begin{figure}[H]
\includegraphics[scale=0.35]{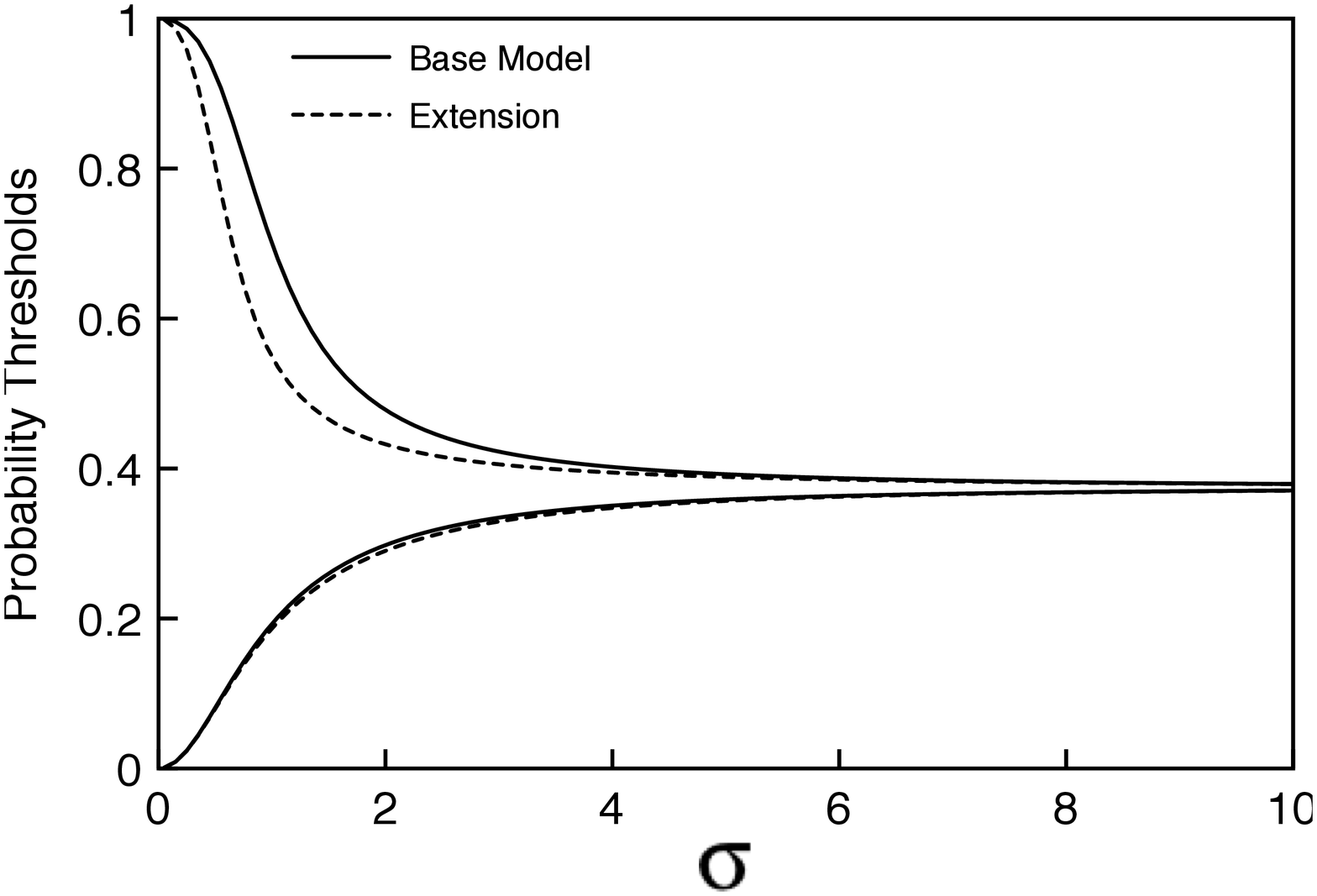}\includegraphics[scale=0.35]{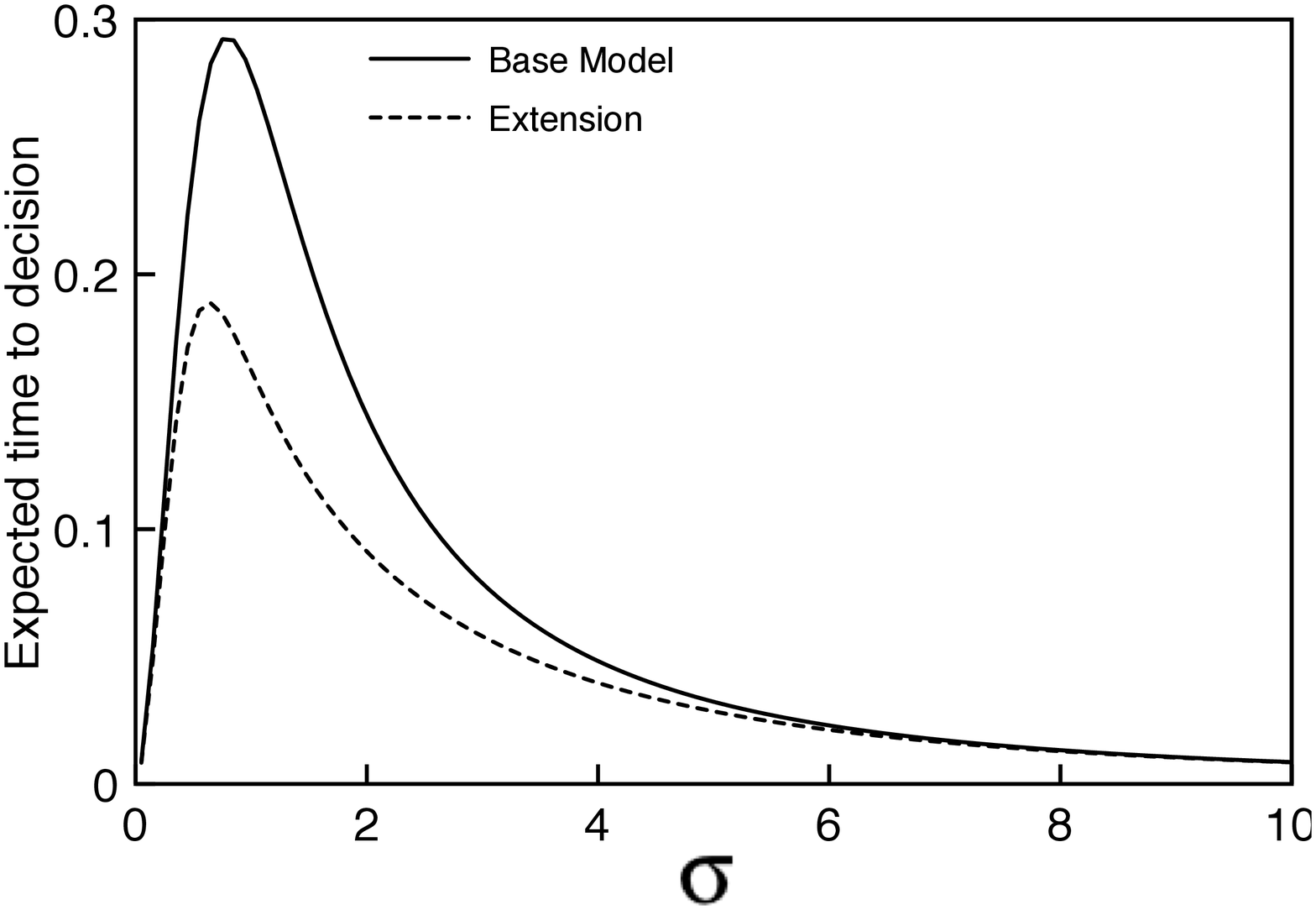}

\caption{}
\label{fig:ex1}
\end{figure}

\begin{figure}[H]
\includegraphics[scale=0.35]{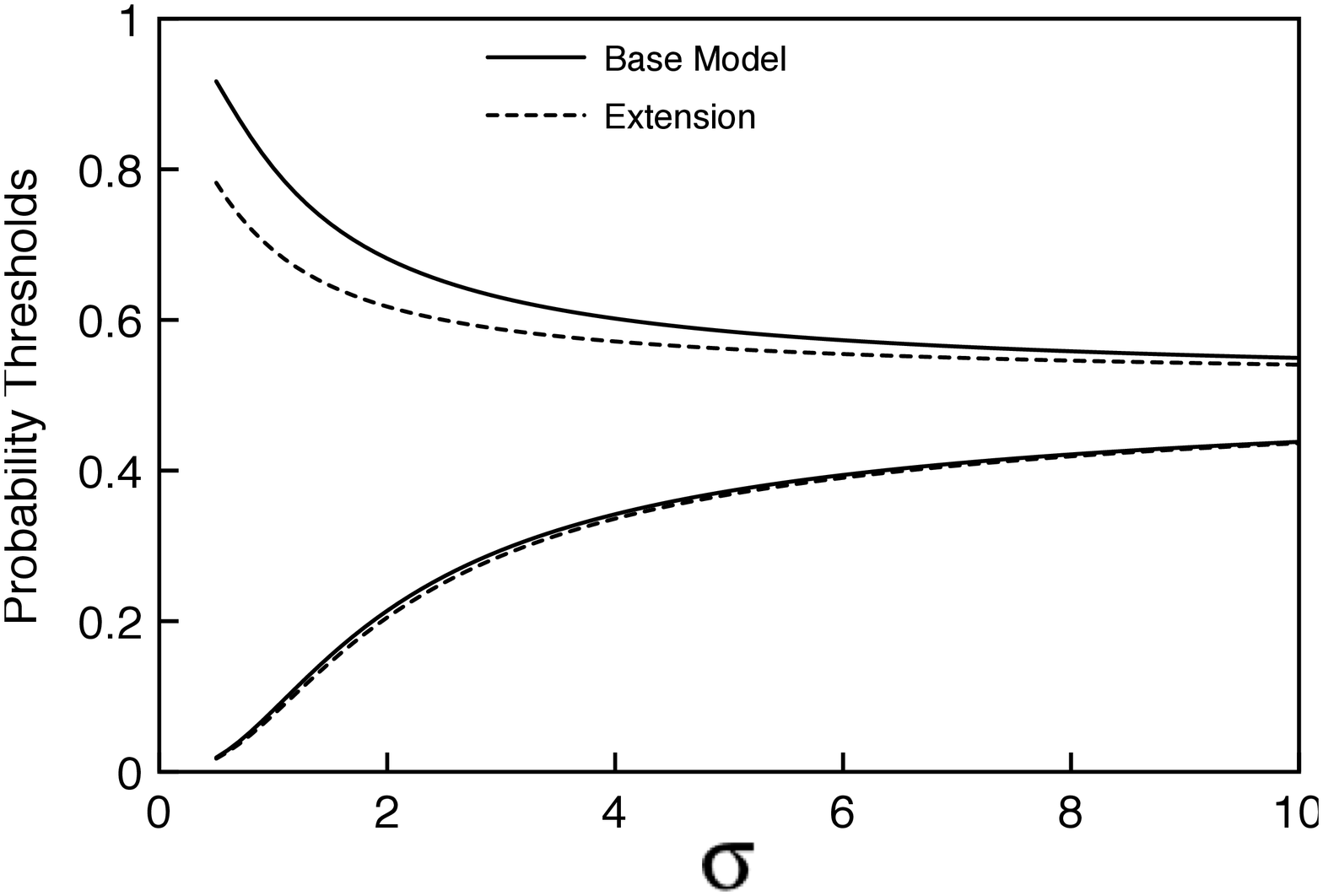}\includegraphics[scale=0.35]{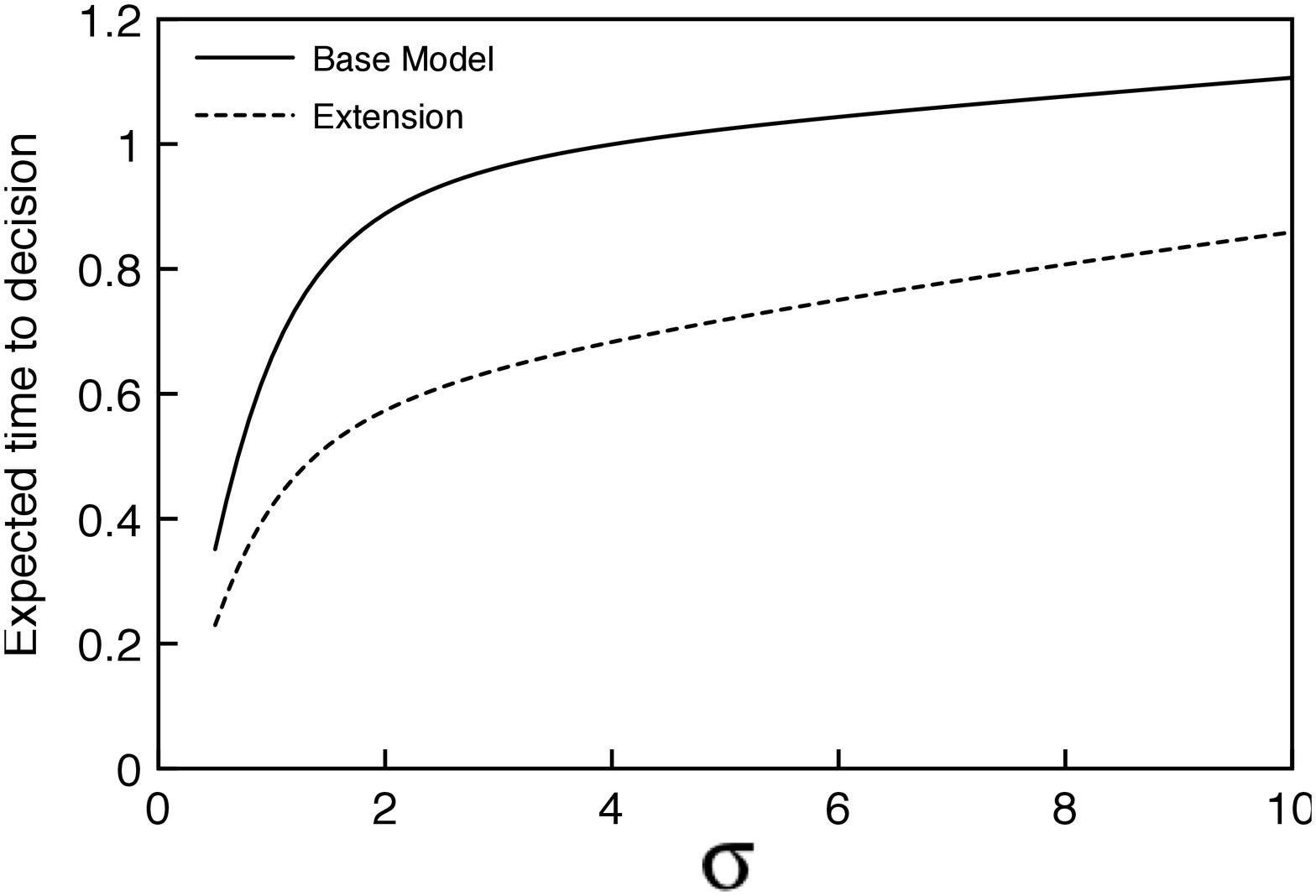}

\caption{}
\label{fig:ex2}
\end{figure}

From Proposition \ref{prop:sigma} and Theorem \ref{thm:Asymptotics},
we can deduce the global behavior of $\underline{\theta}$, $\bar{\theta}$,
and $E^{p}[\tau^{*}]$ as a function of the volatility $\sigma$.
As $\sigma$ increases from zero to very large numbers, $\underline{\theta}$
increases from zero to a limiting value, and $\bar{\theta}$ decreases
from 1 to another limiting value. (See Figs. \ref{fig:ex1} and \ref{fig:ex2}.)
The large-$\sigma$ behavior comports with intuition. When $\sigma$
is very large, arrival of new information is so slow that the firm
should proceed as if $P_{t}$ will never change. In particular, suppose
$g(\hat{p})\ge0$. Then take almost immediate action: exit if $P_{t}<\hat{p}$
and expand if $P_{t}\ge\hat{p}$. Instead, suppose $g(\hat{p})<0.$
Then exit if $P_{t}<\frac{-\ell}{h-\ell}$ and expand if $P_{t}\ge\frac{-\ell^{\dagger}+k\alpha}{h^{\dagger}-\ell^{\dagger}}$;
otherwise, the firm can expect to wait a very long time until it makes
an expansion or exit decision.

For small values of $\sigma$, the expected time to decision $E^{p}[\tau^{*}]$
is very small because it takes a very short time to learn the true
state. The behavior of $E^{p}[\tau^{*}]$ for large $\sigma$ depends
on the initial probability and the sign of $g(\hat{p})$. If $g(\hat{p})\ge0$,
then $g(p)\ge0$ for all $p\in[0,1]$ because $g(\cdot)$ takes its
minimum value at $\hat{p}$. In this case, if $\sigma$ is very large,
due to slow arrival of information, it is not worth spending much
time on the pilot project, and the firm is better off making a quick
expansion or exit decision. Consequently, $E^{p}[\tau^{*}]$ is very
small for large values of $\sigma$. As illustrated in Fig. \ref{fig:ex1}
the expected time-to-decision initially increases in $\sigma$, achieves
a maximum at an intermediate value of $\sigma$, and then approaches
zero for large values of $\sigma$. If $g(\hat{p})<0$, then $g(p)<0$
if and only if $p\in(\frac{-\ell}{h-\ell},\frac{-\ell^{\dagger}+k\alpha}{h^{\dagger}-\ell^{\dagger}})$.
In this case, as shown in Fig. \ref{fig:ex2}, $E^{\hat{p}}[\tau^{*}]$
keeps increasing for large $\sigma$ because it takes a long time
for $P_{t}$ to exit the interval $(\frac{-\ell}{h-\ell},\frac{-\ell^{\dagger}+k\alpha}{h^{\dagger}-\ell^{\dagger}})$
due to very low speed of evolution of $P_{t}$. 

The case $g(\hat{p})<0$ arises, for instance, when $h^{A}/h=\ell^{A}/\ell>0$.
The relation $h^{A}/h=\ell^{A}/\ell>0$ applies to Example 1. If the
production plant is expanded $(n+1)$-fold, then the profit per unit
time is also expanded $(n+1)$-fold: $h^{A}=(n+1)h$ and $\ell^{A}=(n+1)\ell$,
so $h^{A}/\ell^{A}=h/\ell$. 

\section{Exit Option After Expansion \label{sec:Exit-Option}}

From both a legal and a practical perspective, the ability to abandon
a project, before or after an expansion, is always available though
there are costs associated with abandonment. Such an embedded exit
option is known to result in unconventional comparative statics of
the optimal policy with respect to the volatility \citep{Alvarez2003a,Kwon2010}.
Hence, we are motivated to study the effect of the post-expansion
exit option. In this section, we investigate the effect of the post-expansion
exit option upon the optimal policy and the time-to-decision, and
we scrutinize the robustness of our comparative statics results of
Theorem \ref{thm:Asymptotics}. We show that the comparative statics
in the limiting cases of small and large values of $\sigma$ are robust
against the embedded exit option. 

\subsection{The Objective Function}

In this model with an embedded option, the decision-maker's policy
determines two stopping times: $\tau_{A}$ (time to acquire assets)
and $\tau_{E}$ (time to exit). The cumulative profit process is $X_{t}$
given by Eq. (\ref{eq:X-t}) until $\hat{\tau}\equiv\min\{\tau_{A},\tau_{E}\}$.
After acquisition at time $\tau_{A}$, the cumulative profit process
is represented by $X_{t}^{A}$ where
\[
dX_{t}^{A}=\mu^{A}dt+\sigma^{A}dB_{t}\:.
\]
Here, $\sigma^{A}=(n+1)\sigma$ for Example 1 (expansion) and $\sigma^{A}=\sigma$
for Example 2 (acquisition of a dedicated facility). Then, the objective
function is given by
\[
V_{\tau_{A},\tau_{E}}(p_{0})=E^{p_{0}}[\int_{0}^{\hat{\tau}}e^{-\alpha t}dX_{t}+\mathbf{1}_{\{\tau_{A}<\tau_{E}\}}(-ke^{-\alpha\tau_{A}}+\int_{\tau_{A}}^{\tau_{E}}e^{-\alpha t}dX_{t}^{A})]\:.
\]
From the strong Markov property of $\{B_{t}:t\ge0\}$, we find that
\[
V(p_{0})\equiv\sup_{\tau_{A},\tau_{E}}V_{\tau_{A},\tau_{E}}(p_{0})=\sup_{\tau}E^{p_{0}}[\int_{0}^{\tau}e^{-\alpha t}dX_{t}+e^{-\alpha\tau}\max\{V^{A}(P_{\tau})-k,0\}]\:,
\]
where 
\[
V^{A}(x)=\sup_{\tau}E^{x}[\int_{0}^{\tau}e^{-\alpha t}dX_{t}^{A}]\:.
\]
In analogy with Eq. (\ref{eq:V}), we re-express $V(\cdot)$ as
\[
V(x)=\frac{1}{\alpha}[xh+(1-x)\ell]+\sup_{\tau}R_{\tau}(x)
\]
where
\begin{eqnarray*}
R_{\tau}(x) & = & E^{x}[e^{-\alpha\tau}g(P_{\tau})]\:,\\
g(x) & = & \frac{1}{\alpha}\max\{\alpha V^{A}(x)-k\alpha-xh-(1-x)\ell,-xh-(1-x)\ell\}\:.
\end{eqnarray*}

If we assume that $(h^{A}-\ell^{A})/\sigma^{A}=(h-\ell)/\sigma$,
which is the case of Examples 1 and 2, then the form of $V^{A}(\cdot)$
is obtained by \citet{Ryan2003} as follows: 
\[
V^{A}(x)=\left\{ \begin{array}{ll}
\frac{1}{\alpha}[xh^{A}+(1-x)\ell^{A}]+C\cdot x^{\frac{1}{2}(1-\gamma)}(1-x)^{\frac{1}{2}(1+\gamma)} & \mbox{if}\;x>\theta^{A}\\
0 & \mbox{otherwise}
\end{array}\right.\,.
\]
 where 
\begin{eqnarray}
\theta^{A} & = & \frac{-(\gamma-1)\ell^{A}}{(\gamma+1)h^{A}-(\gamma-1)\ell^{A}}\;,\nonumber \\
C & = & \frac{2}{\alpha\sqrt{\gamma^{2}-1}}(\frac{\gamma-1}{\gamma+1})^{\gamma/2}(h^{A})^{\frac{1}{2}(1-\gamma)}(-\ell^{A})^{\frac{1}{2}(1+\gamma)}\:.\label{eq:ext-C}
\end{eqnarray}
For the remainder of this section, we assume that $V^{A}(1)>k$ so
that it is profitable to expand when the posterior probability $p$
is sufficiently high.

\subsection{Optimal Policy}

In this subsection, we obtain the optimal stopping time $\tau^{*}$
that maximizes $R_{\tau}(\cdot)$. We do so by establishing the existence
of $\tau^{*}$ and the necessary conditions for the continuation region. 

In analogy with Sec. \ref{subsec:Policy}, we define $\hat{p}$ as
the point at which $V^{A}(\hat{p})-k=0$, and we focus on the interesting
case of $\hat{p}\in(0,1)$. 

\begin{theorem} \label{thm:ext-policy} Assume $\hat{p}\in(0,1)$.
The optimal policy to maximize $R_{\tau}(\cdot)$ always exists. Moreover,
the continuation region includes a component $(\underline{\theta},\bar{\theta})$
such that $0\le\underline{\theta}<\hat{p}<\bar{\theta}\le1$. The
optimal return function $R^{*}(\cdot)\equiv R_{\tau^{*}}(\cdot)$
in the interval $[\underline{\theta},\bar{\theta}]$ is given by

\begin{equation}
R^{*}(p)=\left\{ \begin{array}{ll}
c_{1}\psi(p)+c_{2}\phi(p) & \mbox{if }\:p\in(\underline{\theta},\bar{\theta})\\
g(p) & \mbox{if }\:p\in\{\underline{\theta},\bar{\theta}\}\:,
\end{array}\right.\label{eq:ext-Rp}
\end{equation}
where $c_{1}$ and $c_{2}$ are positive numbers such that $R^{*}(\cdot)$
is continuously differentiable. \end{theorem} Note that this theorem
is an analog of Theorem \ref{thm:basic-policy}(i) and that Theorem
\ref{thm:basic-policy}(ii) also applies to this model without modification.
Theorem \ref{thm:ext-policy} does not preclude the existence of other
components of the continuation region disconnected from $(\underline{\theta},\bar{\theta})$
because it is technically difficult to preclude them. Even if other
components of the continuation region exist, the interval $[\underline{\theta},\bar{\theta}]$
is the only region where the dichotomous (exit vs expansion) decision
takes place, so we will focus only on the region $[\underline{\theta},\bar{\theta}]$.

\subsection{Comparative Statics \label{subsec:ext-ComparativeStatics}}

In this subsection, we scrutinize the robustness of the comparative
statics results of Sec. \ref{subsec:ComparativeStatics} and investigate
the effect of the exit option after expansion. The post-expansion
exit option has the following effect on $\underline{\theta}$:

\begin{proposition} \label{prop:ext-lower-t}The post-expansion exit
option induces $\underline{\theta}$ to weakly decrease. 

\end{proposition}

If the post-expansion exit option is available, then the expansion
option becomes more attractive to the decision-maker, and it is optimal
to wait longer before making a permanent exit from the pilot project.
In contrast, the effect of the post-expansion exit option on the upper
threshold $\bar{\theta}$ is not mathematically straightforward. Hence,
the effect on the expected time-to-decision is even less clear because
it depends on both $\underline{\theta}$ and $\bar{\theta}$. Although
a general analytical result is unavailable regarding this issue, we
can obtain partial answers in the limits of small and large $\sigma$.
We also obtain the comparative statics of $\underline{\theta}$ and
$\bar{\theta}$ in the same limiting cases.

\begin{lemma} \label{lem:ext-small-s} In the small-$\sigma$ limit,
\begin{eqnarray}
\underline{\theta} & = & \sigma^{2}\frac{2\alpha}{(h-\ell)^{2}}\frac{(-\ell)}{h^{A}-k\alpha}+o(\sigma^{2})\;,\label{eq:ext-lt-small-s}\\
\bar{\theta} & = & 1-\sigma^{2}\frac{2(h^{\dagger}-k\alpha)}{k(h-\ell)^{2}}+o(\sigma^{2})\:.\label{eq:ext-ut-small-s}
\end{eqnarray}

\end{lemma} 

We also note that the upper threshold is $1-\sigma^{2}\frac{2(h^{\dagger}-k\alpha)}{(k-\ell^{A}/\alpha)(h-\ell)^{2}}+o(\sigma^{2})$
if the post-expansion exit option is absent. \footnote{A function $f(\sigma)$ is said to be $o(\sigma^{2})$ if $f(\sigma)/\sigma^{2}\rightarrow0$
in the limit as $\sigma\rightarrow0$. } In contrast, the lower threshold is not affected by the post-expansion
exit option up to $o(\sigma^{2})$. Hence, the following proposition
obtains:

\begin{proposition} \label{prop:ext-small-s}In the small-$\sigma$
limit, the post-expansion exit option causes both $\bar{\theta}$
and the expected time-to-decision to decrease. \end{proposition}

With a post-expansion exit option, the expansion option is more attractive
to the decision-maker, so it makes sense that the threshold for expansion
is lower. The expected time-to-decision is also smaller because the
decision-maker makes an earlier decision to expand the project due
to a higher value of the expansion option. It is also understandable
that the effect of Proposition \ref{prop:ext-lower-t} (smaller value
of $\underline{\theta}$) is a secondary effect on the expected time-to-decision:
the increased value of the expansion option does not make an exit
option too much less attractive to the decision-maker when the posterior
probability $P_{t}$ is low (when the prospect of exercising the expansion
option is low). Numerical examples are displayed in Figs. \ref{fig:ex1}
and \ref{fig:ex2} which show that the extra exit option indeed decreases
$\underline{\theta}$, $\bar{\theta}$, and the expected time-to-decision. 

In the large-$\sigma$ limit, we obtain the following result:

\begin{proposition} \label{prop:ext-large-s} In the large-$\sigma$
limit, the post-expansion exit option does not have any effect on
$\underline{\theta}$, $\bar{\theta}$, or the expected time-to-decision
up to $O(\sigma^{-n})$ for any finite positive $n$. \end{proposition}

The above proposition merely states that there is essentially no effect
of the post-expansion exit option in the large-$\sigma$ limit. This
is consistent with the intuition that the option value from the learning
opportunity is minimal if $\sigma$ is very large due to the slow
arrival of information, and hence the effect of the post-expansion
exit option is also minimal.

Lastly, we are in a position to inspect the comparative statics of
the expected time-to-decision with respect to $\sigma$. The next
theorem follows immediately from Lemma \ref{lem:ext-small-s} and
Proposition \ref{prop:ext-large-s}:

\begin{theorem}[Expected Time-to-Decision]\label{thm:time-p} Theorem
\ref{thm:Asymptotics} (i), (ii), and (iii) holds for the model with
a post-expansion exit option. \end{theorem} Figures \ref{fig:ex1}
and \ref{fig:ex2} also illustrate a numerical example of Theorem
\ref{thm:time-p}. Therefore, we established the robustness of our
result to the additional exit option. 

\section{Reducing the Volatility\label{sec:Discussion}}

In some business situations, the firm has an external source of information
regarding the profitability of the project that effectively reduces
$\sigma$. Therefore, the comparative statics with respect to $\sigma$
can provide a useful prediction of how the additional source of information
impacts the optimal policy and the duration $\tau^{*}$ of the pilot
project. 

Suppose that the firm observes another firm's profit stream from a
similar project subject to the same uncertainty of the state. Assume
that the other firm's cumulative profit process is 
\[
X_{t}^{e}=\mu t+\sigma^{e}B_{t}^{e}\:,
\]
where the volatility $\sigma^{e}$ of the \emph{external} information
process is known to the firm but not necessarily the same as $\sigma$.
Here $B_{t}^{e}$ is an unobservable one-dimensional Brownian motion.
Moreover, assume that $X_{t}^{e}$ is independent of $X_{t}$ except
that they share the same drift $\mu$. 

Define the \emph{reduced} volatility $\sigma^{r}\equiv\sigma\sigma^{e}/\sqrt{\sigma^{2}+(\sigma^{e})^{2}}$,
a new one-dimensional standard Brownian motion $B_{t}^{r}\equiv\sigma^{r}[B_{t}/\sigma+B_{t}^{e}/\sigma^{e}]$,
and a new process $X_{t}^{r}\equiv(\sigma^{r})^{2}[X_{t}/\sigma^{2}+X_{t}^{e}/(\sigma^{e})^{2}]=\mu t+\sigma^{r}B_{t}^{r}$.
The new Brownian motion $B^{r}$ is unobservable, but the process
$X^{r}$ is constructed from observing both $X$ and $X^{e}$, so
$X^{r}$ is an observable process. Let $\{\mathcal{F}_{t}^{r}:t\ge0\}$
be the natural filtration with respect to the process $X^{r}$. Then
we can construct the posterior process $P_{t}^{r}\equiv\mathcal{P}[\mu=h\vert\mathcal{F}_{t}^{r}]$
which is adapted to $\{\mathcal{F}_{t}^{r}:t\ge0\}$. Using Baye's
rule, we can show that the posterior belief is given by 
\begin{equation}
P_{t}^{r}\equiv\mathcal{P}^{p}[\mu=h\vert\mathcal{F}_{t}^{r}]=\frac{p\exp[\frac{h}{(\sigma^{r})^{2}}X_{t}^{r}-\frac{1}{2}(\frac{h}{\sigma^{r}})^{2}t]}{p\exp[\frac{h}{(\sigma^{r})^{2}}X_{t}^{r}-\frac{1}{2}(\frac{h}{\sigma^{r}})^{2}t]+(1-p)\exp[\frac{\ell}{(\sigma^{r})^{2}}X_{t}^{r}-\frac{1}{2}(\frac{\ell}{\sigma^{r}})^{2}t]}\:,\label{eq:P-r}
\end{equation}
 which is clearly constructed from the observable process $X^{r}$.
As before, the posterior process satisfies the stochastic differential
equations 
\[
dP_{t}^{r}=\frac{h-\ell}{\sigma^{r}}P_{t}^{r}\cdot(1-P_{t}^{r})d\tilde{W}_{t}\:,
\]
 where $\tilde{W}_{t}$ is another one-dimensional standard Brownian
motion which is observable. Lastly, the objective function for this
problem is
\begin{equation}
V_{\tau}^{r}(p)=E^{p_{}}[\int_{0}^{\tau}e^{-\alpha t}dX_{t}+e^{-\alpha\tau}r(P_{\tau}^{r})]=\frac{1}{\alpha}[ph+(1-p)\ell]+E^{p}[e^{-\alpha\tau}g(P_{\tau}^{r})]\:,\label{eq:Vr}
\end{equation}
 where $p$ is the initial prior $P_{0}^{r}\equiv\mathcal{P}[\mu=h\vert\mathcal{F}_{0}^{r}]$
and $\tau$ is a stopping time for the filtration $\{\mathcal{F}_{t}^{r}\}$.
The objective function has the same form as Eq. (\ref{eq:V}) except
that the posterior process $P_{t}^{r}$ in Eq. (\ref{eq:Vr}) is constructed
from $X^{r}$ with volatility $\sigma^{r}$ via Eq. (\ref{eq:P-r}).
Therefore, the problem with external information reduces to that of
Sec. \ref{sec:The-Model} with new volatility $\sigma^{r}$ such that
$\sigma^{r}<\min\{\sigma,\sigma^{e}\}$. The only effect of the additional
source of information is to reduce the volatility of the cumulative
profit. 

\section{Conclusions \label{sec:Conclusions}}

Prior to expanding a new project, it behooves the firm to learn more
about the project's profitability. Firms undertaking a pilot project
can learn about profitability and subsequently make expansion and
exit decisions. In this paper, we studied the expansion and exit problem
under incomplete information. Our results show that  there is a non-monotonic
relationship between the expected time-to-decision and  the size of
the continuation region for an irreversible decision: as the volatility
of the cumulative profit increases (so the rate of information arrival
decreases), the optimal continuation region shrinks, but the time-to-decision
does not monotonically decrease. 

Our paper focused on extracting useful insights from a simple and
tractable model. Hence, we have not addressed practical applications
of our findings. In order to model realistic business expansion decisions,
a number of extensions need to be undertaken. For example, our model
assumes that the noise in the cumulative profit is a Wiener process
with constant volatility. Relaxing this assumption would sacrifice
the analytical tractability. More significantly, our model assumes
that the number of possible states of the project is two. For practical
applications, we need to provide an analytical or numerical solution
to problems with a larger number of states. 

As shown by \citet{Decamps2005}, our model can generalize to a problem
in which the number of possible states of the project is larger than
two but finite. {[}See also \citet{Olsen1992} and \citet{Hu1998}
for multidimensional stopping time problems.{]} In such a multi-state
extension, it is rather straightforward to show that much of Sec.
\ref{sec:The-Model} generalizes to multi-state analogs. In particular,
there is a multi-state analog of Theorem \ref{thm:Asymptotics} regarding
the asymptotic behavior of the expected time-to-decision. However,
analytical solutions such as in Eq. (\ref{eq:Optimal-return}) are
not available for a multi-state model. Hence, for practical purposes,
we need to establish numerical procedures to obtain the optimal policy
and the solution. The difficulty is that a numerical procedure would
require continuous differentiability (smooth-pasting condition) of
the optimal solution (\citealt{Muthuraman2008}), but the multi-state
model does not necessarily result in continuously differentiable solutions\footnote{Continuous differentiability over the boundary points requires regularity
of boundary points of the continuation region, but the regularity
is sometimes violated for a degenerate diffusive process such as  a
multi-dimensional generalization of the posterior process $P_{t}$
governed by the dynamics of one-dimensional Brownian motion $\tilde{B}_{t}$.
We thank Renming Song for pointing this out to us.}. An even more difficult problem is a model with a continuous distribution
of states of the project. Finding appropriate numerical solution procedures
for multi-state extensions would constitute a major research program.

\section*{Acknowledgment}

We thank three anonymous referees and the associate editor for their
helpful suggestions which considerably improved our manuscript. 

\bibliographystyle{INFORMS}
\bibliography{IncompleteInfo}

\begin{thebibliography}{34}
\providecommand{\natexlab}[1]{#1}

\bibitem[{Alvarez(2001)}]{Alvarez2001}
Alvarez, L. H.~R. 2001. Reward functionals, salvage values and optimal
  stopping. \emph{Mathematical Methods of Operations Research} \textbf{54}
  315--337.

\bibitem[{Alvarez(2003)}]{Alvarez2003}
Alvarez, L. H.~R. 2003. On the properties of r-excessive mappings for a class
  of diffusions. \emph{Annals of Applied Probability} \textbf{13} 1517--1533.

\bibitem[{Alvarez and Stenbacka(2004)}]{Alvarez2003a}
Alvarez, L. H.~R., R.~Stenbacka. 2004. Optimal risk adoption: a real options
  approach. \emph{Economic Theory} \textbf{23}(1) 123--147.

\bibitem[{Arrow(1962)}]{Arrow1962}
Arrow, K.~J. 1962. The economic implications of learning by doing. \emph{The
  Review of Economic Studies} \textbf{29}(3) 155--173.

\bibitem[{Bergemann and Välimäki(2000)}]{Bergemann2000}
Bergemann, D., J.~Välimäki. 2000. Experimentation in markets. \emph{The Review
  of Economic Studies} \textbf{67}(2) 213--234.

\bibitem[{Bernardo and Chowdhry(2002)}]{Bernardo2002}
Bernardo, A.~E., B.~Chowdhry. 2002. Resources, real options, and corporate
  strategy. \emph{Journal of Financial Economics} \textbf{63}(2) 211--234.

\bibitem[{Bolton and Harris(1999)}]{Bolton1999}
Bolton, P., C.~Harris. 1999. Strategic experimentation. \emph{Econometrica}
  \textbf{67}(2) 349--374.

\bibitem[{Borodin and Salminen(2002)}]{Borodin2002}
Borodin, A.~N., P.~Salminen. 2002. \emph{Handbook of Brownian Motion -- Facts
  and Formulae}. Birkhauser, Basel, 2nd ed.

\bibitem[{Dayanik and Karatzas(2003)}]{Dayanik2003}
Dayanik, S., I.~Karatzas. 2003. On the optimal stopping problem for
  one-dimensional diffusions. \emph{Stochastic Processes and their
  Applications} \textbf{107} 173--212.

\bibitem[{Decamps et~al.(2005)Decamps, Mariotti, and Villeneuve}]{Decamps2005}
Decamps, J.-P., T.~Mariotti, S.~Villeneuve. 2005. Investment timing under
  incomplete information. \emph{Mathematics of Operations Research}
  \textbf{30}(2) 472--500.

\bibitem[{Decamps et~al.(2006)Decamps, Mariotti, and Villeneuve}]{Decamps2006}
Decamps, J.-P., T.~Mariotti, S.~Villeneuve. 2006. Irreversible investment in
  alternative projects. \emph{Economic Theory} \textbf{28}(2) 425--448.

\bibitem[{Dixit(1989)}]{Dixit1989}
Dixit, A. 1989. Entry and exit decisions under uncertainty. \emph{The Journal
  of Political Economy} \textbf{97}(3) 620--638.

\bibitem[{Dixit(1992)}]{Dixit1992}
Dixit, A. 1992. Investment and hysteresis. \emph{The Journal of Economic
  Perspectives} \textbf{6}(1) 107--132.

\bibitem[{Ghemawat and Stander(1998)}]{Ghemawat1998}
Ghemawat, P., H.~Stander. 1998. Nucor at a crossroads. \emph{Harvard Business
  School Case} 9--793--039.

\bibitem[{Hu and Oksendal(1998)}]{Hu1998}
Hu, Y., B.~Oksendal. 1998. Optimal time to invest when the price processes are
  geometric brownian motions. \emph{Finance and Stochastics} \textbf{2}(3)
  295--310.

\bibitem[{Keller and Rady(1999)}]{Keller1999}
Keller, G., S.~Rady. 1999. Optimal experimentation in a changing environment.
  \emph{The Review of Economic Studies} \textbf{66}(3) 475--507.

\bibitem[{Kwon(2010)}]{Kwon2010}
Kwon, H.~D. 2010. {Invest or Exit? Optimal Decisions in the Face of a Declining
  Profit Stream}. \emph{Operations Research} \textbf{58}(3) 638--649.

\bibitem[{Kwon et~al.(2008)Kwon, Lippman, and Tang}]{Kwon2008a}
Kwon, H.~D., S.~A. Lippman, C.~S. Tang. 2008. When to adjust price under
  incomplete information. UCLA Working Paper.

\bibitem[{Lai(2001)}]{Lai2001}
Lai, T.~L. 2001. Sequential analysis: Some classical problems and new
  challenges. \emph{Statistica Sinica} \textbf{11} 303--408.

\bibitem[{Lippman and Rumelt(1992)}]{Lippman1992}
Lippman, S.~A., R.~P. Rumelt. 1992. Demand uncertainty, capital specificity,
  and industry evolution. \emph{Ind Corp Change} \textbf{1}(1) 235--262.

\bibitem[{McCardle(1985)}]{McCardle1985}
McCardle, K.~F. 1985. Information acquisition and the adoption of new
  technology. \emph{Management Science} \textbf{31}(11) 1372--1389.

\bibitem[{Moscarini and Smith(2001)}]{Moscarini2001}
Moscarini, G., L.~Smith. 2001. The optimal level of experimentation.
  \emph{Econometrica} \textbf{69}(6) 1629--1644.

\bibitem[{Muthuraman and Kumar(2008)}]{Muthuraman2008}
Muthuraman, K., S.~Kumar. 2008. Solving free-boundary problems with
  applications in finance. \emph{Foundations and Trends in Stochastic Systems}
  \textbf{1}(4) 259--341.

\bibitem[{Oksendal(2003)}]{Oksendal03}
Oksendal, B. 2003. \emph{Stochastic Differential Equations: An Introduction
  with Applications}. Springer, 6th ed.

\bibitem[{Olsen and Stensland(1992)}]{Olsen1992}
Olsen, T.~E., G.~Stensland. 1992. On optimal timing of investment when cost
  components are additive and follow geometric diffusions. \emph{Journal of
  Economic Dynamics and Control} \textbf{16}(1) 39--51.

\bibitem[{Peskir and Shiryaev(2006)}]{Peskir2006}
Peskir, G., A.~Shiryaev. 2006. \emph{Optimal Stopping and Free-Boundary
  Problems}. Birkhauser Basel.

\bibitem[{Poor and Hadjiliadis(2008)}]{Poor2008}
Poor, H.~V., O.~Hadjiliadis. 2008. \emph{Quickest Detection}. Cambridge
  University Press.

\bibitem[{Ryan and Lippman(2003)}]{Ryan2003}
Ryan, R., S.~A. Lippman. 2003. Optimal exit from a project with noisy returns.
  \emph{Probab. Engrg. Inform. Sci.} \textbf{17}(04) 435--458.

\bibitem[{Shiryaev(1978)}]{Shiryaev1978}
Shiryaev, A. 1978. \emph{Optimal Stopping Rules}. Springer-Verlag, New York.

\bibitem[{Shiryaev(1967)}]{Shiryaev1967}
Shiryaev, A.~N. 1967. Two problems of sequential analysis. \emph{Cybernetics
  and Systems Analysis} \textbf{3}(2) 63--69.

\bibitem[{Ulu and Smith(2009)}]{Ulu2009}
Ulu, C., J.~E. Smith. 2009. {Uncertainty, Information Acquisition, and
  Technology Adoption}. \emph{Operations Research} \textbf{57}(3) 740--752.

\bibitem[{Wald(1945)}]{Wald1945}
Wald, A. 1945. Sequential tests of statistical hypotheses. \emph{Annals of
  Mathematical Statistics} \textbf{16} 117--186.

\bibitem[{Wald(1973)}]{Wald1973}
Wald, A. 1973. \emph{Sequential Analysis}. Dover Publications.

\bibitem[{Wang(2005)}]{Wang2005}
Wang, H. 2005. A sequential entry problem with forced exits. \emph{Mathematics
  of Operations Research} \textbf{30} 501--520.

\end{thebibliography}

\vfill\pagebreak

\vfill\pagebreak
{} \vfill\pagebreak

\section*{Online Appendix}

\noindent \textbf{Proof of Proposition \ref{prop:singe-threshold}:}
Statement (i) is proved by \citet{Ryan2003} and \citet{Kwon2008a}.
Statement (ii) can be proved similarly. \eproof

\noindent \textbf{Proof of Theorem \ref{thm:basic-policy}:} (i) Consider
the following test function:
\[
f(p)=\left\{ \begin{array}{ll}
c_{1}\psi(p)+c_{2}\phi(p) & \mbox{if }\:p\in(\underline{\theta},\bar{\theta})\\
g(p) & \mbox{otherwise}\:.
\end{array}\right.
\]
 By Theorem 10.4.1 of \citet{Oksendal03}, the test function $f(p)$
is the optimal return if the following conditions are satisfied: (a)
$f(\cdot)$ is continuously differentiable on $[0,1]$, (b) $f(p)\ge g(p)$
for all $p\in[0,1]$, (c) $f(\cdot)$ is twice continuously differentiable
except at $\{\underline{\theta},\bar{\theta}\}$, (d) the second-order
derivatives of $f(\cdot)$ are finite near $\underline{\theta}$ and
$\bar{\theta}$, and (e) $\mathcal{A}f(p)\le0$ for $p\in[0,1]\backslash\{\bar{\theta},\underline{\theta}\}$.
(Note that $\mathcal{A}f(p)$ is not defined at $p=\bar{\theta}$
or $\underline{\theta}$.)  Conditions (c) and (d) can be readily
verified from the form of $f(\cdot)$. It remains to show that there
exist $c_{1},c_{2},\underline{\theta},$ and $\bar{\theta}$ which
satisfy (a), (b), and (e).

(a) The condition that $f(\cdot)$ is continuously differentiable
can be expressed as follows:
\begin{eqnarray}
c_{1}\psi(\underline{\theta})+c_{2}\phi(\underline{\theta}) & = & -\alpha^{-1}[\underline{\theta}h+(1-\underline{\theta})\ell]\,,\label{eq:bdry1}\\
c_{1}\psi(\bar{\theta})+c_{2}\phi(\bar{\theta}) & = & \alpha^{-1}[\bar{\theta}h^{\dagger}+(1-\bar{\theta})\ell^{\dagger}]-k\,,\label{eq:bdry2}\\
c_{1}\psi'(\underline{\theta})+c_{2}\phi'(\underline{\theta}) & = & -\alpha^{-1}(h-\ell)\:,\label{eq:pasting1}\\
c_{1}\psi'(\bar{\theta})+c_{2}\phi'(\bar{\theta}) & = & \alpha^{-1}(h^{\dagger}-\ell^{\dagger})\:.\label{eq:pasting2}
\end{eqnarray}
Our strategy is to prove that there always exists a solution $\underline{\theta},\bar{\theta},c_{1}$,
and $c_{2}$ that satisfies Eqs. (\ref{eq:bdry1})-(\ref{eq:pasting2})
and the inequality $\underline{\theta}<\bar{\theta}$. 

For algebraic convenience, we introduce $\underline{q}\equiv\underline{\theta}/(1-\underline{\theta})$,
$\bar{q}\equiv\bar{\theta}/(1-\bar{\theta})$, and $\beta\equiv\bar{q}/\underline{q}$,
where $0<\underline{q}\le\bar{q}<\infty$ and $\beta\in[1,\infty)$.
After eliminating $c_{1}$ and $c_{2}$ from Eqs. (\ref{eq:bdry1})-(\ref{eq:pasting2}),
we obtain a pair of equations for $\underline{q}$ and $\beta$ as
follows:
\begin{eqnarray}
\underline{q} & = & Q_{1}(\beta)\equiv(\frac{\gamma-1}{\gamma+1})\frac{-\ell-(\ell^{\dagger}-k\alpha)\beta^{-\frac{1}{2}(\gamma+1)}}{h+(h^{\dagger}-k\alpha)\beta^{\frac{1}{2}(1-\gamma)}}\ \label{eq:Q1-b}\\
 & = & Q_{2}(\beta)\equiv(\frac{\gamma+1}{\gamma-1})\frac{-\ell-(\ell^{\dagger}-k\alpha)\beta^{\frac{1}{2}(\gamma-1)}}{h+(h^{\dagger}-k\alpha)\beta^{\frac{1}{2}(\gamma+1)}}\;.\label{eq:Q2-b}
\end{eqnarray}
We note that $Q_{1}(1)<Q_{2}(1)$ but $\lim_{\beta\rightarrow\infty}Q_{1}(\beta)=(\frac{-\ell}{h})(\frac{\gamma-1}{\gamma+1})>0$
and $\lim_{\beta\rightarrow\infty}Q_{2}(\beta)=0$, so there is always
a solution $\beta\in(1,\infty)$ to $Q_{1}(\beta)=Q_{2}(\beta)$;
once $\beta$ is determined, $\underline{q}$ is determined by $\underline{q}=Q_{1}(\beta)$. 

(b) Note that the fundamental solutions $\psi(\cdot)$ and $\phi(\cdot)$
are monotonically increasing and decreasing respectively, and they
are both convex. Hence, if $c_{1}$ and $c_{2}$ have opposite signs,
then the test function $f(p)=c_{1}\psi(p)+c_{2}\phi(p)$ in the interval
$(\underline{\theta},\bar{\theta})$ must be either monotonically
increasing or monotonically decreasing; this contradicts Eqs. (\ref{eq:pasting1})
and (\ref{eq:pasting2}) because $h-\ell>0$ and $h^{\dagger}-\ell^{\dagger}>0$.
Thus, $c_{1}$ and $c_{2}$ are either both positive or both negative,
so $f(\cdot)$ must be either convex or concave because the fundamental
solutions $\psi(\cdot)$ and $\phi(\cdot)$ are both convex. However,
$f(\cdot)$ cannot be concave because $f'(\bar{\theta})>0>f'(\underline{\theta})$
by Eqs. (\ref{eq:pasting1}) and (\ref{eq:pasting2}). We conclude
that $f(\cdot)$ must be convex. It follows that $c_{1}$ and $c_{2}$
are both positive. 

From the fact that $g(\cdot)$ is linear except at $\hat{p}$ and
that the first derivative of $g(\cdot)$ discontinuously increases
at $\hat{p}$, it follows that $\underline{\theta}<\hat{p}<\bar{\theta}$
because otherwise the first derivative of $f(\cdot)$ is constant
over the interval $(\underline{\theta},\bar{\theta})$. Moreover,
because $f''(p)>g''(p)$ for all $p\in(\underline{\theta},\bar{\theta})\backslash\{\hat{p}\}$,
we have the inequality $f(p)>g(p)$ for $p\in(\underline{\theta},\bar{\theta})$. 

(e) Note that $f(\cdot)$ is a strictly positive function. Hence,
$g(p)$ is positive for $p\not\in(\underline{\theta},\bar{\theta})$.
Because $g(\cdot)$ is linear except at $\hat{p}\in(\underline{\theta},\bar{\theta})$
(piecewise linear), the inequality $\mathcal{A}g(p)=-\alpha g(p)<0$
is satisfied for $p\not\in(\underline{\theta},\bar{\theta})$. Hence,
$\mathcal{A}f(p)=0$ for $p\in(\underline{\theta},\bar{\theta})$
and $\mathcal{A}f(p)=\mathcal{A}g(p)<0$ for $p\not\in(\underline{\theta},\bar{\theta})$.
This concludes the proof of part (i).

(ii) In order to compute $E^{p}[\tau^{*}]$, we define $\mathcal{L}$,
the infinitesimal generator of $P_{t}$, by
\[
\mathcal{L}\equiv\frac{1}{2}\frac{(h-\ell)^{2}}{\sigma^{2}}p^{2}(1-p)^{2}\frac{\partial^{2}}{\partial p^{2}}
\]
and notice that the probability that $P_{\tau}$ hits the upper threshold
first is given by $\mathcal{P}(P_{\tau^{*}}=\bar{\theta})=(p-\underline{\theta})/(\bar{\theta}-\underline{\theta})$
from II.4 and II.9 of \citet{Borodin2002}. Using the fact that $\mathcal{L}T(p)=1$
and Dynkin's formula (\citealt{Oksendal03}), we obtain 
\begin{eqnarray*}
E^{p}[\tau^{*}] & = & E^{p}[\int_{0}^{\tau^{*}}ds]=E^{p}[\int_{0}^{\tau^{*}}\mathcal{L}T(P_{t})dt]\\
 & = & E^{p}[T(P_{\tau^{*}})]-T(p)\:,
\end{eqnarray*}
which reduces to Eq. (\ref{eq:E-tau}).\eproof

\noindent \textbf{Proof of Proposition \ref{prop:sigma}:} First,
note that the volatility $\frac{h-\ell}{\sigma}P_{t}(1-P_{t})$ of
the process $\mathbf{P}=\{P_{t}:t\ge0\}$ in Eq. (\ref{eq:P-t-SDE})
is inversely proportional to $\sigma$. Second, note that $R^{*}(p)$
is convex in $p$. Now we can directly apply Theorem 4 of \citet{Alvarez2003}
and conclude that $R^{*}(p)$ is non-decreasing in the volatility
of $P_{t}$ and hence non-increasing in $\sigma$. Because the continuation
region $(\underline{\theta},\bar{\theta})$ is defined as $\{p:R^{*}(p)>g(p)\}$,
it follows that $\bar{\theta}$ decreases in $\sigma$ and that $\underline{\theta}$
increases in $\sigma$. \eproof

\noindent \textbf{Proof of Theorem \ref{thm:Asymptotics}:} (i) In
the small-$\sigma$ limit, we expand Eqs. (\ref{eq:Q1-b}) and (\ref{eq:Q2-b})
in powers of $\sigma$ and study the leading-order terms to obtain
the following asymptotic solution to $Q_{1}(\beta)=Q_{2}(\beta)$
where $\beta\equiv\bar{\theta}(1-\underline{\theta})/[\underline{\theta}(1-\bar{\theta})]$
:
\begin{eqnarray*}
\beta & = & \frac{(\ell^{A}-k\alpha)(h^{A}-k\alpha)}{\ell(h^{\dagger}-k\alpha)}\cdot\frac{(h-\ell)^{4}}{4\sigma^{4}\alpha^{2}}+o(\sigma^{-4})\:,\\
\underline{\theta} & = & \frac{-\ell}{h^{A}-k\alpha}\cdot\frac{2\sigma^{2}\alpha}{(h-\ell)^{2}}+o(\sigma^{2})\:,\\
\bar{\theta} & = & 1-\frac{(h^{\dagger}-k\alpha)}{(-\ell^{A}+k\alpha)}\cdot\frac{2\sigma^{2}\alpha}{(h-\ell)}+o(\sigma^{2})\:.
\end{eqnarray*}
 Moreover, from $T(p)=\frac{2\sigma^{2}}{(h-\ell)^{2}}(2p-1)\ln\frac{x}{1-x}$,
we have $T(\underline{\theta})=-\frac{2\sigma^{2}}{(h-\ell)^{2}}\ln\sigma^{2}+O(\sigma^{2})$
and $T(\bar{\theta})=\frac{2\sigma^{2}}{(h-\ell)^{2}}\ln\frac{1}{\sigma^{2}}+O(\sigma^{2})$,
so $E^{p}[\tau^{*}]\rightarrow0$ as $\sigma\rightarrow0$.

(ii) and (iii): From 
\[
g(\hat{p})=\frac{1}{\alpha}\frac{h(\ell^{\dagger}-k\alpha)-\ell(h^{\dagger}-k\alpha)}{h^{A}-\ell^{A}}
\]
 and $h^{A}-\ell^{A}>0$, the condition $g(\hat{p})\ge0$ reduces
to $h(\ell^{\dagger}-k\alpha)-\ell(h^{\dagger}-k\alpha)\ge0$ while
$g(\hat{p})<0$ reduces to $h(\ell^{\dagger}-k\alpha)-\ell(h^{\dagger}-k\alpha)<0$.
From the fact that $\underline{\theta}$ increases and $\bar{\theta}$
decreases in $\sigma$, we find that two alternative cases are possible
in the $\sigma\rightarrow\infty$ limit : $\beta\searrow1$ or $\beta\searrow\beta_{C}>1$. 

Assume that $\beta\searrow\beta_{C}>1$. Solving Eqs. (\ref{eq:Q1-b})
and (\ref{eq:Q2-b}) in the large-$\sigma$ limits, we find 
\[
\beta_{C}=\frac{h(\ell^{\dagger}-k\alpha)}{\ell(h^{\dagger}-k\alpha)}
\]
so that $\beta\searrow\beta_{C}>1$ if and only if $h(\ell^{\dagger}-k\alpha)-\ell(h^{\dagger}-k\alpha)<0$
or $g(\hat{p})<0$. In this case, from Eq. (\ref{eq:Q1-b}), we obtain
the limiting behaviors $\underline{\theta}\nearrow-\ell/(h-\ell)$
and $\bar{\theta}\searrow-(\ell^{\dagger}-k\alpha)/(h^{\dagger}-\ell^{\dagger})$
as $\sigma\rightarrow\infty$. Moreover, from Eq. (\ref{eq:E-tau}),
it can be shown that $E^{p}[\tau^{*}]\rightarrow\infty$ if $p\in(\frac{-\ell}{h-\ell},\frac{-\ell^{\dagger}+k\alpha}{h^{\dagger}-\ell^{\dagger}})$.
On the other hand, if $p\not\in[\frac{-\ell}{h-\ell},\frac{-\ell^{\dagger}+k\alpha}{h^{\dagger}-\ell^{\dagger}}]$,
then $E^{p}[\tau^{*}]\rightarrow0$ because $p$ eventually belongs
to the stopping region (the complement of the continuation region)
as $\sigma\rightarrow\infty$. Lastly, from Eq. (\ref{eq:E-tau})
and the fact that $\beta=\beta_{C}+O(\sigma^{-1})$, $\underline{\theta}=\frac{-\ell}{h-\ell}+O(\sigma^{-1})$,
and $\bar{\theta}=\frac{-\ell^{\dagger}+k\alpha}{h^{\dagger}-\ell^{\dagger}}+O(\sigma^{-1})$,
we obtain $E^{p}[\tau^{*}]=O(\sigma)$ for $p=\frac{-\ell}{h-\ell}$
or $\frac{-\ell^{\dagger}+k\alpha}{h^{\dagger}-\ell^{\dagger}}$. 

If $h(\ell^{\dagger}-k\alpha)-\ell(h^{\dagger}-k\alpha)\ge0$ or $g(\hat{p})\ge0$,
solving Eqs. (\ref{eq:Q1-b}) and (\ref{eq:Q2-b}), we obtain 
\[
\beta=1+\frac{(h^{A}-k\alpha)(\ell^{A}-k\alpha)}{\ell(h^{\dagger}-k\alpha)-h(\ell^{\dagger}-k\alpha)}\cdot\frac{(h-\ell)^{2}}{2\alpha^{2}\sigma^{2}}+O(\sigma^{-4})\:,
\]
 $\underline{\theta}=\hat{p}+O(\sigma^{-2})$, and $\bar{\theta}=\hat{p}+O(\sigma^{-2})$
such that $\underline{\theta}\nearrow\hat{p}$, $\bar{\theta}\searrow\hat{p}$
as $\sigma\rightarrow\infty$. Moreover, from Eq. (\ref{eq:E-tau}),
$E^{\hat{p}}[\tau^{*}]=O(\sigma^{-2})$. \eproof

\noindent \textbf{Proof of Theorem \ref{thm:ext-policy}:} We show
that (a) the optimal policy exists, (b) there is a component of the
continuation region containing $\hat{p}$, and (c) the optimal return
function is a continuously differentiable function satisfying Eq.
(\ref{eq:ext-Rp}). 

(a) We closely follow the proof of Proposition 7.1 in \citet{Decamps2005}.
Upon inspection of the objective function $R_{\tau}(\cdot)$, it is
clear that the optimal policy, if it exists, should be stationary
because the discounted reward function $e^{-\alpha t}g(p)$ is homogeneous
in time. Our objective is to see if there is $\tau^{*}$ such that
\[
\bar{R}(p)\equiv\sup_{\tau}R_{\tau}(p)=R^{*}(p)\:.
\]
 Note that $(t,P_{t})$ is a Feller process and that $e^{-\alpha t}g(p)$
is a continuous and bounded function of $(t,p)$ for $t>0$ and $p\in\Delta$;
hence, $\bar{R}(p)$ is bounded and lower semicontinuous (\citealt{Peskir2006},
p. 49). It follows that the set $C=\{p:\bar{R}(p)>g(p)\}$ is an open
set, and the exit time $\tau_{C}=\inf\{t\ge0:P_{t}\not\in C\}$ is
a well-defined stopping time. From \citet{Decamps2005} p. 497, the
stopped process $\{e^{-\alpha(t\wedge\tau_{C})}\bar{R}(P_{t\wedge\tau_{C}}):t\ge0\}$
is a martingale. Hence, for any positive integer $M$, 
\begin{eqnarray*}
\bar{R}(p) & = & E^{p}[e^{-\alpha(M\wedge\tau_{C})}\bar{R}(P_{M\wedge\tau_{C}})]\\
 & = & E^{p}[e^{-\alpha\tau_{C}}g(P_{\tau_{C}})\mathbf{1}_{\{\tau_{C}\le M\}}]+E^{p}[e^{-\alpha M}\bar{R}(P_{M})\mathbf{1}_{\{\tau_{C}>M\}}]\:.
\end{eqnarray*}
The supremum $\bar{R}(p)$ is always non-negative because the return
from never stopping is $R_{\infty}(p)=0$ and $R_{\infty}(p)\le\bar{R}(p)$
by the definition of $\bar{R}(p)$. By the monotone convergence theorem,
we have 
\[
\lim_{M\rightarrow\infty}E^{p}[e^{-\alpha\tau_{C}}g(P_{\tau_{C}})\mathbf{1}_{\{\tau_{C}\le M\}}]=E^{p}[e^{-\alpha\tau_{C}}g(P_{\tau_{C}})]\:.
\]
 Also 
\[
\lim_{M\rightarrow\infty}E^{p}[e^{-\alpha M}\bar{R}(P_{M})\mathbf{1}_{\{\tau_{C}>M\}}]=0\ 
\]
 because $\bar{R}(\cdot)$ is bounded. Thus, $\bar{R}(p)=E^{p}[e^{-\alpha\tau_{C}}g(P_{\tau_{C}})]=R^{*}(p)$,
so $\tau^{*}=\tau_{C}$: an optimal policy exists.

(b) Because the continuation region is an open set, it suffices to
show that the indifference point $\hat{p}$ belongs to the continuation
region, which was achieved by \citet{Decamps2006} Proposition 2.2.
Hence, there is an interval $(\underline{\theta},\bar{\theta})$ containing
$\hat{p}$ which is a component of the continuation region.

(c) Note that $g(\cdot)$ is a bounded and continuous function, and
$\mathbf{P}$ is a Markov process. Hence, we can apply Theorem 3.15
of \citet{Shiryaev1978}, p. 157 to prove $\mathcal{A}R^{*}(p)=0$
for $p\in C$ and $R^{*}(p)=g(p)$ for $p\notin C$. Thus, for $p\in C$,
$R^{*}(p)$ is given by Eq. (\ref{eq:ext-Rp}). Moreover, because
$\mathbf{P}$ is a diffusive process without any jump, $R^{*}(\cdot)$
is continuous (\citealt{Peskir2006}, p. 148) at $\underline{\theta}$
and $\bar{\theta}$. Lastly, because of the regularity of the diffusive
process $\mathbf{P}$, the first derivative of $R^{*}(\cdot)$ is
continuous at $\underline{\theta}$ and $\bar{\theta}$ (\citealt{Peskir2006},
p. 150). \eproof

\noindent \textbf{Proof of Proposition \ref{prop:ext-lower-t}}: Let
$R_{0}^{*}(\cdot)$ and $(\underline{\theta}_{0},\bar{\theta}_{0})$
denote the optimal return function and the continuation region \emph{without}
the post-expansion exit option. {[}$R_{0}^{*}(\cdot)$ and $(\underline{\theta}_{0},\bar{\theta}_{0})$
are obtained in Sec. \ref{sec:The-Model}.{]} Similarly, let $g_{0}(\cdot)$
denote the gain function defined in Eq. (\ref{eq:g-x}). Note that
$(\underline{\theta}_{0},\bar{\theta}_{0})=\{p\in(0,1):R_{0}^{*}(p)>g_{0}(p)\}$
and $(\underline{\theta},\bar{\theta})=\{p\in(0,1):R^{*}(p)>g(p)\}$.
Because $g(p)\ge g_{0}(p)$ for all $p$, we have $R^{*}(p)\ge R_{0}^{*}(p)$.
Moreover, for sufficiently small $p$, $g(p)=g_{0}(p)$. 

Suppose that $\underline{\theta}_{0}<\underline{\theta}$. Then $g(\underline{\theta})=g_{0}(\underline{\theta})<R_{0}^{*}(\underline{\theta})$.
However, we have $g(\underline{\theta})=R^{*}(\underline{\theta})\ge R_{0}^{*}(\underline{\theta})$,
a contradiction. We conclude that $\underline{\theta}_{0}\ge\underline{\theta}$.
\eproof 

\noindent \textbf{Proof of Lemma \ref{lem:ext-small-s}}: The equations
for $\underline{\theta}$ and $\bar{\theta}$ can be obtained from
the condition that $R^{*}(\cdot)$ is continuously differentiable
at $\underline{\theta}$ and $\bar{\theta}$:
\begin{eqnarray*}
c_{1}\psi(\underline{\theta})+c_{2}\phi(\underline{\theta}) & = & -\alpha^{-1}[\underline{\theta}h+(1-\underline{\theta})\ell]\,,\\
c_{1}\psi(\bar{\theta})+c_{2}\phi(\bar{\theta}) & = & V^{A}(\bar{\theta})-k-\alpha^{-1}[\bar{\theta}h+(1-\bar{\theta})\ell]\,,\\
c_{1}\psi'(\underline{\theta})+c_{2}\phi'(\underline{\theta}) & = & -\alpha^{-1}(h-\ell)\:,\\
c_{1}\psi'(\bar{\theta})+c_{2}\phi'(\bar{\theta}) & = & V^{A\prime}(\bar{\theta})-\alpha^{-1}(h-\ell)\:.
\end{eqnarray*}
After eliminating $c_{1}$ and $c_{2}$, using the form of $V^{A}(p)$
when $V^{A}(p)>k$, we obtain the following equations in terms of
$\underline{q}\equiv\underline{\theta}/(1-\underline{\theta})$, $\bar{q}\equiv\bar{\theta}/(1-\bar{\theta})$,
and $\beta\equiv\bar{q}/\underline{q}$:
\begin{eqnarray}
\underline{q} & = & (\frac{\gamma-1}{\gamma+1})\frac{-\ell-(\ell^{\dagger}-k\alpha)\beta^{-\frac{1}{2}(\gamma+1)}}{h+(h^{\dagger}-k\alpha)\beta^{\frac{1}{2}(1-\gamma)}}\ \label{eq:ext-Q1}\\
 & = & (\frac{\gamma+1}{\gamma-1})\frac{-\ell-(\ell^{\dagger}-k\alpha)\beta^{\frac{1}{2}(\gamma-1)}}{h+(h^{\dagger}-k\alpha)\beta^{\frac{1}{2}(\gamma+1)}}-\frac{2\alpha C\gamma}{\gamma-1}\frac{\underline{q}^{-\frac{1}{2}(\gamma-1)}}{h+(h^{\dagger}-k\alpha)\beta^{\frac{1}{2}(\gamma+1)}}\label{eq:ext-Q2}
\end{eqnarray}
From Eqs. (\ref{eq:ext-Q1}), (\ref{eq:ext-Q2}), and (\ref{eq:ext-C}),
we obtain Eqs. (\ref{eq:ext-lt-small-s}) and (\ref{eq:ext-ut-small-s}).
\eproof 

\noindent \textbf{Proof of Proposition \ref{prop:ext-small-s}}: From
$\ell^{A}<0$, we have $1-\sigma^{2}\frac{2(h^{\dagger}-k\alpha)}{k(h-\ell)^{2}}<1-\sigma^{2}\frac{2(h^{\dagger}-k\alpha)}{(k-\ell^{A}/\alpha)(h-\ell)^{2}}$,
so $\bar{\theta}$ is smaller with the post-expansion exit option.
It follows that the expected time-to-decision is also smaller up to
$o(\sigma^{2})$. \eproof 

\noindent \textbf{Proof of Proposition \ref{prop:ext-large-s}}: Let
$\beta_{\infty}$ denote $\lim_{\sigma\rightarrow\infty}\beta$. Then
we consider two cases: $\beta_{\infty}>1$ and $\beta_{\infty}=1$.

(i) $\beta_{\infty}>1$: From Eq. (\ref{eq:ext-Q1}), the second term
of Eq. (\ref{eq:ext-Q2}) is
\[
-4\frac{\sqrt{h^{A}\ell^{A}\ell/(h\beta_{\infty})}}{e\gamma(h^{\dagger}-k\alpha)}(\frac{\ell^{A}h}{\ell h^{A}\beta_{\infty}})^{\gamma/2}(1+o(1))\:.
\]
For Example 1, $h^{A}/\vert\ell^{A}\vert=h/\vert\ell\vert$, so $\frac{\ell^{A}h}{\ell h^{A}\beta_{\infty}}<1$.
For Example 2, $h^{A}>h$ and $\vert\ell\vert>\vert\ell^{A}\vert$,
so again $\frac{\ell^{A}h}{\ell h^{A}\beta_{\infty}}<1$. It follows
that $(\frac{\ell^{A}h}{\ell h^{A}\beta_{\infty}})^{\gamma/2}$ converges
to zero faster than $\sigma^{-n}$ for any positive $n$. Thus, Eqs.
(\ref{eq:ext-Q1}) and (\ref{eq:ext-Q2}) coincide with (\ref{eq:Q1-b})
and (\ref{eq:Q2-b}) in the large-$\sigma$ limit.

(ii) $\beta_{\infty}=1$: The only possible solution for $\beta$
is of the form $1+O(\sigma^{-2})$. (There is no term of order $O(\sigma^{-1})$.)
The second term of Eq. (\ref{eq:ext-Q2}) is
\[
-\frac{4}{e\gamma}\sqrt{\frac{h^{A}\ell^{A}(\ell^{A}-k\alpha)}{(h^{A}-k\alpha)^{3}}}[\frac{(h^{A}-k\alpha)(-\ell^{A})}{h^{A}(-\ell^{A}+k\alpha)}]^{\gamma/2}(1+o(1))\:.
\]
(The factor of $e^{-1}$ comes from the factor $\lim_{\sigma\rightarrow}(\frac{\gamma-1}{\gamma+1})^{\gamma/2}$
in the definition of $C$.)  From the fact that $\frac{(h^{A}-k\alpha)(-\ell^{A})}{h^{A}(-\ell^{A}+k\alpha)}<1$,
this term converges to zero faster than $\sigma^{-n}$ for any positive
$n$. Thus, Eqs. (\ref{eq:ext-Q1}) and (\ref{eq:ext-Q2}) coincide
with (\ref{eq:Q1-b}) and (\ref{eq:Q2-b}) in the large-$\sigma$
limit. \eproof 
\end{document}